\documentclass[reqno,10pt]{amsart}
\usepackage[foot]{amsaddr}

\usepackage{amssymb,amsmath,amsthm,amsfonts}
\usepackage{mathrsfs,dsfont,comment,mathscinet,mathtools}

\usepackage{subdepth}
\usepackage[T1]{fontenc}
\usepackage{cite}

\usepackage[utf8]{inputenc}
\usepackage[normalem]{ulem}
\usepackage[left=2.5cm,right=2.5cm,top=2.5cm,bottom=2.5cm]{geometry}

\usepackage{randomwalk}
\usepackage{pgfplots}
\usepackage{tkz-fct}

\usepackage{graphicx,tikz,color}
\usetikzlibrary{quotes,angles}

\usepackage[format=hang,labelfont=bf]{caption}

\usepackage{enumerate,esint}

\usepackage{ wasysym }

\usepackage{epstopdf}
\usepackage[colorlinks=true, pdfstartview=FitV, linkcolor=blue, 
            citecolor=blue, urlcolor=blue]{hyperref}

\allowdisplaybreaks
\numberwithin{equation}{section}
\theoremstyle{plain}
\newtheorem{theorem}{Theorem}[section]
\newtheorem{proposition}[theorem]{Proposition}
\newtheorem{lemma}[theorem]{Lemma}

\theoremstyle{definition}
\newtheorem{definition}[theorem]{Definition}

\newtheorem*{theorem*}{Theorem}

\makeatletter
\def\th@plain{%
  \thm@notefont{}
  \itshape 
}
\def\th@definition{%
  \thm@notefont{}
  \normalfont 
}
\makeatother

\let\oldproofname=\proofname
\renewcommand{\proofname}{\bf{\oldproofname}}

\definecolor{mblue}{HTML}{13439b}

\newcommand{\R}{\mathbb R}
\newcommand\N{\mathbb N}

\newcommand{\So}{\mathbb{S}^1}
\newcommand{\Rt}{\mathbb{R}^2}
\newcommand{\Rtt}{\mathbb{R}^{2\times 2}}

\renewcommand{\d}{\mathrm{d}}
\renewcommand{\H}{\mathscr{H}}

\newcommand{\restr}[1]{|_{#1}}

\makeatletter
\newcommand{\subsetsim}{\mathrel{\mathpalette\subset@sim\relax}}
\newcommand{\subset@sim}[2]{%
	\vbox{\offinterlineskip\m@th
		\ialign{\hfil$#1##$\hfil\cr
			\sim\cr\subset\cr
		}%
	}%
}

\usepackage{accents}
\newcommand*{\dt}[1]{%
	\accentset{\mbox{\large\bfseries .}}{#1}}

\def\Xint#1{\mathchoice
{\XXint\displaystyle\textstyle{#1}}%
{\XXint\textstyle\scriptstyle{#1}}%
{\XXint\scriptstyle\scriptscriptstyle{#1}}%
{\XXint\scriptscriptstyle\scriptscriptstyle{#1}}%
\!\int}
\def\XXint#1#2#3{{\setbox0=\hbox{$#1{#2#3}{\int}$ }
\vcenter{\hbox{$#2#3$ }}\kern-.6\wd0}}

\def\dashint{\Xint-}

\makeatletter
\newcommand{\subalign}[1]{%
	\vcenter{%
		\Let@ \restore@math@cr \default@tag
		\baselineskip\fontdimen10 \scriptfont\tw@
		\advance\baselineskip\fontdimen12 \scriptfont\tw@
		\lineskip\thr@@\fontdimen8 \scriptfont\thr@@
		\lineskiplimit\lineskip
		\ialign{\hfil$\m@th\scriptstyle##$&$\m@th\scriptstyle{}##$\hfil\crcr
			#1\crcr
		}%
	}%
}
\makeatother

\definecolor{bblue}{HTML}{3C3C9F}

\newcommand{\cof}{\mathrm{cof}}

\newcommand{\supp}{\mathrm{supp}}
\renewcommand{\deg}{\mathrm{deg}}
\newcommand{\Det}{\mathrm{Det}}
\newcommand{\per}{\mathrm{Per}}

\newcommand{\imt}{\mathrm{im}_{\rm T}}

\renewcommand{\div}{\mathrm{div}}
\newcommand{\dist}{\mathrm{dist}}
\newcommand\wk{\rightharpoonup}

\newcommand{\wks}{\overset{\ast}{\rightharpoonup}}

\newcommand*\closure[1]{\overline{#1}}

\newcommand{\leb}{\mathscr{L}^2}
\newcommand{\haus}{\mathscr{H}^1}

\def\Xint#1{\mathchoice
{\XXint\displaystyle\textstyle{#1}}%
{\XXint\textstyle\scriptstyle{#1}}%
{\XXint\scriptstyle\scriptscriptstyle{#1}}%
{\XXint\scriptscriptstyle\scriptscriptstyle{#1}}%
\!\int}
\def\XXint#1#2#3{{\setbox0=\hbox{$#1{#2#3}{\int}$ }
\vcenter{\hbox{$#2#3$ }}\kern-.6\wd0}}

\def\dashint{\Xint-}

\newcommand\EEE{\color{black}}
\newcommand\RRR{\color{red}}

\title[Quasistatic growth of cavities and cracks]{Quasistatic growth of cavities and cracks in the plane}
\author[M. Bresciani]{Marco Bresciani${}^{*}$}
\address{*\,Department Mathematik, Friedrich-Alexander-Universit\"{a}t Erlangen-N\"{u}rnberg, Cauerstrasse 11, 91058 Erlangen (DE)}
\email{marco.bresciani@fau.de}
\author[M. Friedrich]{Manuel Friedrich${}^{*}$}
\email{manuel.friedrich@fau.de}

\begin{document}

\setlength\parindent{0pt}

\vskip .2truecm
\begin{abstract}
  We propose a model  for the  quasistatic growth of cavities and cracks in two-dimensional nonlinear elasticity. Cavities and cracks are modeled as discrete and compact subsets of a  planar  domain, respectively, and deformations are defined only outside of cracks.  The  model accounts for the irreversibility of both processes of cavitation and fracture and it allows for the coalescence of cavities into cracks.  Our main result shows the  existence of quasistatic evolutions in the case of a finite number of cavities, under an a priori bound on the number of connected components of the cracks.     
\end{abstract}
\maketitle


\section{Introduction}

\subsection{Motivation}
Numerous experiments conducted on both soft elastomers \cite{gent.lindley,gent.park,poulain.etal} and ductile metals \cite{petrinic1,petrinic2} have shown how fracture  is generally preceded by the formation of small voids in the interior of solids. In view of this observation,  fracture initiation can be explained as the outcome of  growth and  coalescence of  such voids.  The   failure mechanism   of void formation in response to external loading \EEE is usually referred to as cavitation and  
  allows materials to evade unsustainable states of stress by undergoing singular deformations. In the past, cavitation has been interpreted as an elastic instability determined by the abrupt expansion of preexisting defects in materials \cite{ball.cavitation,gent.lindley}. More recently,  various investigations have objected this interpretation  on  the basis of empirical evidences (see, e.g., \cite{kumar.lopez-pamies,lefevre} and the references therein). These arguments have led to the conviction that cavitation cannot be completely resolved just on purely elastic grounds, but it inherently involves the irreversible creation of new surface by  deformation whose growth is governed by Griffith's criterion \cite{griffith}.    

In agreement with this  line of thought, a variational model in the context of nonlinear elasticity allowing for both cavitation and fracture has been proposed in \cite{henao.moracorral.invertibility}. The model ultimately  rests on the minimization of an energy functional consisting of three terms: an elastic term depending on the deformation gradient, a term accounting for jump discontinuities of the deformation, and a surface term measuring the area of the new surface created by the deformation. The latter comprises the perimeter of the cavities opened by the deformation and  the area of  the fracture surface  in the deformed configuration, but also the invisible surface occurring when pieces of surfaces created at different points in the reference configuration are put in contact  \cite{henao.moracorral.fracture}.    It has been shown that the boundedness of the surface energy  has  consequences in connection to the weak continuity of Jacobian determinants \cite{henao.moracorral.invertibility} and the regularity of inverse deformations \cite{henao.moracorral.fracture,henao.moracorral.regularity}.  

The existence theory in \cite{henao.moracorral.invertibility} provides minimizers of the energy functional  under physical growth conditions, thus ensuring the well posedness of the model in the static setting. The present work originates  from  the intention of extending these results to the evolutionary setting. Specifically, we are  interested   in cavitation and fracture processes driven by slow loadings and in the resulting quasistatic evolutions.

\subsection{Related literature}
\label{subsec:lit}
In formulating our model, we  draw  upon several works. Without aiming for completeness, we briefly mention here only those which are particularly relevant for the present study.

A huge body of literature is devoted to the study of quasistatic evolutions for brittle fracture  (see  \cite{bourdin.etal} for a broad overview).  Following  the modeling approach pioneered in \cite{francfort.marigo}, which   is based on Griffith's theory of crack propagation  \cite{griffith},  evolutions   are determined by the competition between a bulk elastic energy and a term proportional to the crack surface in the reference configuration. Most importantly, the irreversibility of the fracture process is taken into account, namely, the crack set has to grow in time.  As a consequence,   quasistatic evolutions have to comply with the   principle of  {unilateral stability}  (here,  ``unilateral'' is related to irreversibility)  and  additionally with a  power balance.   Their definition corresponds to the notion of  energetic solutions  in the realm of  rate-independent processes \cite[Subsection 4.2.4.1]{mielke.roubicek}.

In this framework, two different main categories of works   can be distinguished, namely   those adopting   \emph{strong} and \emph{weak} formulations. In the strong formulation, \EEE  cracks are modeled as compact subsets of the closure of the reference domain and deformations are defined as Sobolev maps on the cracked domain. As discussed later, technical limitations  confine these results to the two-dimensional setting  and pose some restrictions on the topology of the cracks. In this setting,  antiplanar  \cite{dalmaso.toader, DM-Toa2}  and linear elasticity \cite{acanfora.ponsiglione,chambolle} have been addressed.   We mention also the paper \cite{francfort.giacomini.lopez} which      focuses  on the process of healing within the context of crack propagation,  motivated  by the experimental evidence  that internal cracks arising from   cavitation processes  can actually heal.  In particular, in contrast  to  the previous papers, this work does not enforce any irreversibility, so that cracks might grow or shrink in time.

  In the   {weak formulation},   cracks are described as rectifiable subsets of the domain, while deformations are modeled as (generalized) special maps of bounded variation (denoted by $(G)SBV$, see, e.g., \cite{ambrosio.fusco.pallara}) on the same domain. 
 This setting is more flexible, for  it comes along without any topological assumptions on the crack geometry and, in most cases, models can be treated in arbitrary space dimensions. In particular,  existence results in this setting have been achieved in  \cite{francfort.larsen} for antiplanar elasticity,  in  \cite{friedrich.solombrino} for   two-dimensional   linear elasticity, and  in  \cite{dalmaso.francfort.toader,dalmaso.lazzaroni,dalmaso.lazzaroni2,lazzaroni} for nonlinear elasticity. 

The literature on evolutionary problems for cavitation is scarce. To the best of our knowledge, excluding studies where the location or the geometry of cavities are restricted,  
the quasistatic growth of cavities  has been solely investigated  in \cite{mora.corral}.
The setting of this paper resembles the weak formulation  for the quasistatic growth of  cracks described above. The energy depends on  a discrete set  of flaw points   in the reference domain and a Sobolev deformation  defined on the entire reference domain. The flaw points correspond to a set of points where  deformations \emph{might} open cavities.
The energy functional comprises an elastic contribution and a surface energy measuring the perimeter of the cavities in the deformed configuration. Additionally, the energy includes a dissipative term proportional to the number of flaw points standing for an initiation energy. In particular, deformations with finite energy can open only a finite number of cavities.
 Irreversibility is enforced by requiring that the  set of flaw points grows in time, while  
the set of cavitation points  of deformations, namely, those points at which deformations \emph{actually} open cavities,  is contained in the set of flaw points along the evolution. Also here, quasistatic evolutions must fulfill an unilateral stability condition together with a power balance.

\subsection{The proposed model}
\label{subsec:model} As we already mentioned, the surface energy proposed in \cite{henao.moracorral.invertibility} accounts for both cavitation and fracture.   Yet, in order to extend the  evolutionary model for cavitation \cite{mora.corral}  to the setting of fracture, a rigorous  definition of flaw points and cavities is required. Such a rigorous definition   is possible    for Sobolev deformations with suitable integrability \cite{henao.moracorral.lusin}  and   makes use of topological degree theory \cite{fonseca.gangbo}. However, since the degree is only available for maps enjoying some kind of continuity, the definition cannot be extended to deformations in $SBV$. Also, the naive approach of defining cavitation points  outside of the closure of the jump set as done for Sobolev maps is doomed to fail  since,  in general, jump sets of  $SBV$ functions  might be very irregular, e.g.,    they might in principle even be a dense subset of the domain.  For this reason, it seems currently out of reach to work entirely in a weak formulation for both cracks and cavities. \EEE

Our proposed model for investigating the growth of cavities and fractures features a formulation which sits   between a weak formulation and a strong one.   The    setting is two-dimensional and the   energy functional depends on three variables: a discrete set of flaw points in the closure of the reference domain, a compact subset of the closure of the reference domain modeling the crack  (with an a priori bounded  number  of connected components),   and a Sobolev deformation defined outside of cracks.
The energy  accounts for the elastic energy, the area of the perimeter of the cavities in the deformed configuration, the cardinality of the set of flaw points, and the area of the crack. Therefore, as in \cite{mora.corral}, only deformations opening a finite number of cavities can have finite energy.  
 Despite its limitations,   the model advanced here still allows for the coalescence of cavities into cracks  \cite{henao.moracorral.xu,henao.moracorral.xu2}.     To our view,   the possibility of capturing this  phenomenon  represents  an important and desirable  feature of the model.

In accordance with the literature mentioned in Subsection \ref{subsec:lit}, except for \cite{francfort.giacomini.lopez}, the guiding principles behind our notion of quasistatic evolution are three:  irreversibility, global unilateral stability, and power balance. Irreversibility requires both  the set of flaw points  and  the crack to grow in time. Moreover, as in \cite{mora.corral}, the set of cavitation points of admissible deformations needs to be contained in the set of flaw points at each time. The unilateral stability requires the minimality with respect to competitors with larger or equal set of flaw points and cracks. Eventually, the power, i.e., the time derivative of the energy along the evolution, needs to balance the work of external and displacement loadings.

\subsection{Main result}
Like in \cite{dalmaso.lazzaroni2}, we focus on the pure Neumann problem and we study evolutions  driven by time-dependent  body forces. Moreover, we prescribe an a priori $L^\infty$ bound on the deformations in order to be able to work with (globally) Sobolev maps. 
The main result of the paper is given in Theorem \ref{thm:existence}  and shows the existence of quasistatic evolutions  for  the model outlined in Subsection \ref{subsec:model}. 

The general scheme of the proof of Theorem \ref{thm:existence} is standard and works by time discretization. To solve the incremental minimization problem, we apply the direct method of the calculus of variations. In this regard,  the topology on the class of admissible states plays a crucial role. For  the  sets of flaw points and cracks, we consider the metric topology induced by the Hausdorff distance, while we use the weak topology for deformations and their gradients. The choice of the Hausdorff topology for the cracks mainly motivates the restriction to the two-dimensional case as well as the a priori bound on the number of their connected components. Indeed, it is known that the lower semicontinuity of the  Hausdorff measure   with respect to the Hausdorff convergence of closed sets  can be  guaranteed just in this setting, see, e.g., \cite[Example 4.4.18]{ambrosio.tilli} and \cite[p.~216]{chambolle}. 

A technical aspect that is worth to mention concerns condition (INV). Introduced in \cite{mueller.spector}, this is a topological condition on deformations which  is closely related to almost everywhere injectivity. From the modeling point of view, condition (INV) excludes unphysical deformations and hence constitutes a desirable restriction to pose on admissible deformations. From the mathematical point of view, it enables us to resort to the representation formulas for  the surface energy  and the distributional determinant in \cite[Theorem 4.6]{henao.moracorral.lusin} which we use to obtain controls on the perimeter of the cavities and the position of cavitation points as in \cite[Lemma 3.11]{bresciani.friedrich.moracorral}. In this work, admissible deformations are required to satisfy condition (INV) in the cracked domain. The classical proof of the stability of condition (INV) with respect to the weak convergence of deformations \cite[Lemma~3.3]{mueller.spector} is easily adapted to the variant on cracked domains, thanks to a well-known localization property of Hausdorff convergence. The same property allows us to  deduce the lower semicontinuity of the surface energy for deformations on cracked domains from the original case on intact domains \cite[Theorem 3]{henao.moracorral.invertibility}.


The most technical part of the proof concerns the transfer of cavities and cracks given in Theorem \ref{thm:transfer}. This result represents the counterpart of the famous jump transfer lemma for the quasistatic growth of brittle fracture \cite{dalmaso.francfort.toader,dalmaso.lazzaroni,dalmaso.lazzaroni2,francfort.giacomini.lopez,francfort.larsen,friedrich.solombrino,lazzaroni}, which is used to show the unilateral stability of the candidate solution, in the present setting. In the language of rate-independent processes, this theorem ensures the existence of mutual recovery sequences \cite{mielke.roubicek}.  The proof of Theorem \ref{thm:transfer} demands the explicit construction of a sequence of admissible states and hinges on a careful combination of techniques borrowed from the literature. Specifically, the definition of the sets of flaw points as well as the transfer of the cavities follows the arguments in \cite{mora.corral}. The construction of the cracks exploits a technical lemma from \cite{dalmaso.toader} related to the Hausdorff convergence of compact sets under unilateral constraints, while the crack transfer combines the local-separation property of cracks in \cite{francfort.giacomini.lopez} with an adaptation of the stretching argument in \cite{dalmaso.lazzaroni,dalmaso.lazzaroni2,lazzaroni}. In particular, through  the overall construction, special attention is needed to ensure that the deformations constructed  satisfy condition (INV).   


\subsection{Possible extensions}  In this paper, we did not impose any boundary conditions on the deformations. This choice was made to avoid further technicalities, and  allows us to   focus  on the interplay between cavitation and fracture. However, we do not see any substantial obstacle in including Dirichlet boundary conditions into the model as the necessary techniques are already available in the literature. In particular,  
the approach in \cite{francfort.mielke,dalmaso.lazzaroni,lazzaroni} seems to be   viable also for our problem. 

Even though acceptable within the context of nonlinear elasticity, the a priori $L^\infty$ bound on admissible deformations can be avoided. Following \cite{dalmaso.toader}, this extension can be achieved by working with homogeneous Sobolev spaces (also known as Deny-Lions spaces \cite{deny.lions}), that is, spaces of locally Sobolev maps whose gradient is globally integrable.  

Eventually, having specific applications to  soft organic solids in mind \cite{dougan.etal,jansen.etal}, it would be interesting to study a variant of the present model allowing for the healing of both cavities and cracks. The analysis in \cite{francfort.giacomini.lopez} paves the way in this direction. Suitable modifications of our arguments along the lines of this work seem feasible, but    are   left to future investigations.

\subsection{Structure of the paper}
The paper contains four sections, the first of which constitutes this introduction. In Section \ref{sec:setting-main-result}, we describe the setting and we state our main result. After recalling all the necessary  results from the literature in Section \ref{sec:preliminaries}, the proof of our main result is given in Section \ref{sec:proof}.

\EEE

\section{Setting  and main result}
\label{sec:setting-main-result}

In this section, we describe the setting of the problem and we state our main result.  Some technical definitions needed for the formulation of the model are deferred to Section \ref{sec:preliminaries}.  
  
\subsection{Notation}
Given $a,b\in\R$, we set $a \wedge b\coloneqq \min\{a,b\}$.  
We  work in two dimensions and we indicate the canonical basis of $\Rt$ by $(\boldsymbol{e}_1,\boldsymbol{e}_2)$.  We use the symbols $\Rtt$ for the set of two-by-two matrices and $\Rtt_+$ for those with positive determinant.   
The identity map over $\Rt$ is denoted by $\boldsymbol{id}$, while $\boldsymbol{0}\in\Rt$ is the origin. Also, $\boldsymbol{I}$ and $\boldsymbol{O}$ indicate the identity matrix in $\Rtt$ and the null tensor in $\R^{2\times 2\times 2}$, respectively.

Closure and boundary of a set $A\subset \Rt$ are denoted by $\closure{A}$ and $\partial A$, respectively. The symbol $B(\boldsymbol{x},r)$ is used for the open ball with center $\boldsymbol{x}\in\Rt$ and radius $r>0$. Then, we set $\closure{B}(\boldsymbol{x},r)\coloneqq \closure{B(\boldsymbol{x},r)}$ and $S(\boldsymbol{x},r)\coloneqq \partial B(\boldsymbol{x},r)$. In particular, $\So\coloneqq S(\boldsymbol{0},1)$. \EEE  Given $\boldsymbol{\nu}\in \So$,  we define
\begin{equation*}
	Q_{\boldsymbol{\nu}}(\boldsymbol{x},r)\coloneqq \left\{ \boldsymbol{z}\in\Rt\colon \, \:|(\boldsymbol{z}-\boldsymbol{x})\cdot \boldsymbol{\nu}|<\frac{r}{2},\:|(\boldsymbol{z}-\boldsymbol{x})\cdot \boldsymbol{\nu}^\perp|<\frac{r}{2} \right\},
\end{equation*}
where $\boldsymbol{\nu}^\perp\coloneqq (\boldsymbol{\nu}\cdot \boldsymbol{e}_2,-\boldsymbol{\nu}\cdot \boldsymbol{e}_1)^\top$. That is, $Q_{\boldsymbol{\nu}}(\boldsymbol{x},r)$ is the open cube centered at $\boldsymbol{x}$ with side length $r$  with  two sides parallel to $\boldsymbol{\nu}$. 

We use the symbols $\leb$ and $\H^\alpha$ for the Lebesgue measure and the $\alpha$-dimensional Hausdorff measure on $\Rt$, where $\alpha>0$. We employ standard notation for Lebesgue, Sobolev, and Bochner spaces. When $A\subset \Rt$ is a measurable set, we indicate with $L^1_+(A)$ the class of functions in $L^1(A)$ that are positive almost everywhere.
For a Banach space $X$ and $T>0$, we denote by $AC([0,T];X)$ the space of absolutely continuous functions from $[0,T]$ to $X$. When $X=\R$, we simply write $AC([0,T])$. Given $\boldsymbol{f}\in AC([0,T];X)$, we use the symbol $\dt{\boldsymbol{f}}$ for its derivative which is a function in $L^1(0,T;X)$. Given a measurable set $A\subset \Rt$, we denote its perimeter by $\per(A)$, see, e.g.,  \cite[Section 3.5]{ambrosio.fusco.pallara}.

\subsection{Setting}
Let $\Omega\subset \R^2$ be a bounded Lipschitz domain.     In particular, $\Omega$ is connected.  We denote by $\mathcal{C}(\closure{\Omega})$ and $\mathcal{K}(\closure{\Omega})$ the collections of all finite and compact subsets of $\closure{\Omega}$, respectively.   Given $m\in \N$, we define
\begin{align*}
	\mathcal{K}_m(\closure{\Omega})&\coloneqq \left\{ K \subset \closure{\Omega}:\:\text{$K$ compact with at most $m$ connected components},\:\H^1(K)<+\infty  \right\}.
\end{align*}  
Then, we set 
\begin{equation*}
	\mathcal{X}_m(\closure{\Omega})\coloneqq \mathcal{C}(\closure{\Omega}) \times \mathcal{K}_m(\closure{\Omega}).
\end{equation*}
Elements of $\mathcal{C}(\closure{\Omega})$ are the sets of points  where  deformations can open cavities, while  $\mathcal{K}_m(\closure{\Omega})$ constitutes the class of admissible cracks.  
Let $p  \in (1,2)$ and $M>0$. For every $(C,K)\in\mathcal{X}_m(\closure{\Omega})$, we  define the corresponding class of admissible deformations  as  
\begin{equation*}
	\begin{split}
	\mathcal{Y}_{p,M}(C,K)\coloneqq \Big\{ \boldsymbol{y}\in W^{1,p}(\Omega \setminus K;\R^2)\colon \quad  \|\boldsymbol{y}\|_{L^\infty(\Omega\setminus K;\R^2)}\leq M,\quad   \text{$\det D \boldsymbol{y}\in L^1_+(\Omega \setminus K)$}, \quad \mathcal{S}(\boldsymbol{y})<+\infty&,\\
	  \text{$\boldsymbol{y}$ satisfies  (INV) in $\Omega\setminus K$},\quad C_{\boldsymbol{y}} \subset \RRR C \EEE &\Big \}.
	\end{split}
\end{equation*}
In the previous  definition,  the bound determined by $M$ enforces a confinement condition and the assumption on the sign of the Jacobian determinant restricts the analysis to orientation-preserving deformations. The functional $\mathcal{S}$ is introduced in Definition~\ref{def:surface-energy}, and the value of
$\mathcal{S}(\boldsymbol{y})$ corresponds  to  the  length  of the new surface created by  the cavities opened by  $\boldsymbol{y}$. Condition (INV) is a topological condition which excludes  unphysical deformations. While this condition is usually imposed for  deformations  defined over the intact domain $\Omega$, here it is formulated for  deformations   on the cracked domain $\Omega \setminus K$, see Definition \ref{def:INV} below.   (As   in the classical case, deformations  satisfying (INV) are almost everywhere injective,  see the comment after that definition.)    Eventually,  the set 
$C_{\boldsymbol{y}} \subset \Omega \setminus K$ denotes the points where $\boldsymbol{y}$ opens a cavity in the sense of Definition \ref{def:cavitation-points} below.  

The resulting class of admissible states is given by
\begin{equation*}
	\mathcal{A}_{p,M}^m\coloneqq \big\{ (C,K,\boldsymbol{y}):\:(C,K)\in \mathcal{X}_m(\closure{\Omega}),\:\boldsymbol{y}\in\mathcal{Y}_{p,M}(C,K) \big\}.
\end{equation*}
Let $T>0$ denote the time horizon. The energy functional $\mathcal{F}\colon [0,T]\times \mathcal{A}_{p,M}^m \to \R \cup \{ +\infty\}$ is defined as
\begin{equation}
	\label{eq:F}
	\begin{split}
		\mathcal{F}(t,C,K,\boldsymbol{y})\coloneqq 
		\mathcal{W}(\boldsymbol{y})+\mathcal{P}(C,K,\boldsymbol{y})+\H^0(C)+\H^1(K)-\mathcal{L}(t,\boldsymbol{y}),
	\end{split}
\end{equation}   
where we set
\begin{align}
	\label{eq:WPL}
	\mathcal{W}(\boldsymbol{y})\coloneqq \int_{\Omega \setminus K} W(D\boldsymbol{y})\,\d \boldsymbol{x}, \quad \mathcal{P}(C,K,\boldsymbol{y})\coloneqq\sum_{\boldsymbol{a}\in C\cap (\Omega \setminus K)} \per \left( \imt(\boldsymbol{y},\boldsymbol{a})\right), \quad \mathcal{L}(t,\boldsymbol{y})\coloneqq\int_{\Omega \setminus K} \boldsymbol{f}(t)\cdot \boldsymbol{y}\,\d\boldsymbol{x}.
\end{align}
  The first term in \eqref{eq:F} stands for the elastic energy and depends on the density   $W\colon \Rtt \to [0,+\infty]$.  The  second term is the surface energy that accounts for the stretching of the boundary of the cavities, with the symbol $\imt(\boldsymbol{y},\boldsymbol{a})$ denoting the topological image of a point, as given in Definition \ref{def:cavitation-points}  below.  The third and the fourth term in \eqref{eq:F} are the contributions related to nucleation of cavities and the propagation of cracks, respectively. Eventually, the last term represents the work of applied loads described by  the function   $\boldsymbol{f}\in AC([0,T];L^1(\Omega;\R^2))$.  

On the density $W$   in \eqref{eq:WPL},  we make the following assumptions:
\begin{enumerate}[(W1)]
	\item \emph{Finiteness and regularity}: $W(\boldsymbol{F})=+\infty$ for all $\boldsymbol{F}\in \Rtt$ with $\det \boldsymbol{F}\leq 0 $ and $W$ is continuously differentiable in $\Rtt_+$.
	\item \emph{Growth}: There exist a constant $c_W>0$ and a Borel function $\gamma\colon (0,+\infty)\to [0,+\infty]$  satisfying
	\begin{equation}
		\label{eq:gamma}
		\lim_{h \to 0^+} \gamma(h)= \lim_{h \to +\infty} \frac{\gamma(h)}{h}=+\infty
	\end{equation}
	such that 
	\begin{equation*}
		W(\boldsymbol{F})\geq c_W|\boldsymbol{F}|^p+\gamma(\det \boldsymbol{F}) \quad \text{for all $\boldsymbol{F}\in \Rtt_+$.}
	\end{equation*}
	\item \emph{Polyconvexity}: $W$ is  polyconvex,   that is, there  exists a convex function  $\Phi\colon \Rtt \times  (0,+\infty)  \to [0,+\infty]$  such that 
	\begin{equation*}
		 W(\boldsymbol{F})=\Phi(\boldsymbol{F},\det \boldsymbol{F})   \quad \text{for all $\boldsymbol{F}\in \Rtt_+$.}
	\end{equation*}
	\item \emph{Stress control}: There exist two constants $a_W>0$ and $b_W\geq 0$ such that
	\begin{equation*}
		|\boldsymbol{F}^\top DW(\boldsymbol{F})|\leq a_W(W(\boldsymbol{F})+b_W) \quad \text{for all $\boldsymbol{F}\in \Rtt_+$.}
	\end{equation*}
\end{enumerate}
All these assumption are physical and by now standard,   see, e.g., \cite{ball.op,dalmaso.lazzaroni,francfort.mielke}. We mention that the estimate in (W4) entails an analogous control for the Kirchhoff stress tensor $\boldsymbol{F}\mapsto DW(\boldsymbol{F})\boldsymbol{F}^\top$, under the assumption of frame indifference on  $W$.

\subsection{Main result}
 The notion of quasistatic evolution given below combines similar notions defined in the modeling of cavitation in elastic solids \cite{mora.corral} and fracture of brittle materials in the so-called strong formulation  \cite{ bourdin.etal,    dalmaso.toader,francfort.giacomini.lopez}.  The definition falls into the class of rate-independent evolutions \cite{mielke.roubicek} and accounts for the irreversibility of both processes of cavitation and fracture. \EEE

\begin{definition}[Irreversible quasistatic evolution]
	\label{def:irreversible-qs}
A function $t\mapsto (C(t),K(t),\boldsymbol{y}(t))$ from $[0,T]$ to $\mathcal{A}_{p,M}^m$ is termed an \emph{irreversible quasistatic evolution} if the following  holds: 
\begin{enumerate}[(i)]
	\item \emph{Irreversibility:} for every $t_1,t_2 \in [0,T]$ with $t_1\leq t_2$, we have $C(t_1)\subset C(t_2)$ and $K(t_1)\subset K(t_2)$;
	\item \emph{Unilateral global stability:} for every $t\in [0,T]$, we have
	\begin{equation*}
		\mathcal{F}(t,C(t),K(t),\boldsymbol{y}(t))\leq \mathcal{F}(t,\widetilde{C},\widetilde{K},\widetilde{\boldsymbol{y}})\quad \text{for all $(\widetilde{C},\widetilde{K},\widetilde{\boldsymbol{y}})\in\mathcal{A}_{p,M}^m$ with $\widetilde{C} \supset C(t)$ and $\widetilde{K}\supset K(t)$;}
	\end{equation*}
	\item \emph{Power balance:} the function $t\mapsto  \mathcal{F}(t,C(t),K(t),\boldsymbol{y}(t))$ belongs to $AC([0,T])$ and satisfies
	\begin{equation*}
		\mathcal{F}(t,C(t),K(t),\boldsymbol{y}(t))=\mathcal{F}(0,C(0),K(0),\boldsymbol{y}(0))-\int_0^t \int_{\Omega \setminus K(s)} \dt{\boldsymbol{f}}(s)\cdot \boldsymbol{y}(s)\,\d\boldsymbol{x}\,\d s \quad \text{for all $t\in [0,T]$.}
	\end{equation*}
	\EEE 
\end{enumerate}
\end{definition}

The main result of this paper  shows  the existence of quasistatic  evolutions  as defined above.

\begin{theorem}[Existence of irreversible quasistatic evolutions]
	\label{thm:existence}
Let $(C_0,K_0,\boldsymbol{y}_0)\EEE\in\mathcal{A}_{p,M}^m$ be such that 
\begin{equation}
	\label{eq:initial-datum}
	\mathcal{F}(0,  C_0,K_0,\boldsymbol{y}_0\EEE )\leq \mathcal{F}(0,\widetilde{C},\widetilde{K},\widetilde{\boldsymbol{y}}) \quad \text{for all $(\widetilde{C},\widetilde{K},\widetilde{\boldsymbol{y}})\in\mathcal{A}_{p,M}^m$ with $\widetilde{C}\supset  C_0\EEE $ and $\widetilde{K}\supset  K_0\EEE$.}
\end{equation}	
Then, there exists an irreversible quasistatic evolution $t\mapsto (C(t),K(t),\boldsymbol{y}(t))$ as in Definition \ref{def:irreversible-qs}   satisfying  the initial condition $(C(0),K(0),\boldsymbol{y}(0))=( C_0,K_0,\boldsymbol{y}_0\EEE )$.
\end{theorem}

 Solutions $t\mapsto (C(t),K(t),\boldsymbol{y}(t))$  are not expected to be regular. 
In addition to the monotonicity of $t\mapsto C(t)$ and $t\mapsto K(t)$ coming from  Definition \ref{def:irreversible-qs}(i), by using measurable selections,  it is possible to establish the existence of an evolution for which the map $t\mapsto (\boldsymbol{y}(t),D\boldsymbol{y}  (t)  )$ from $[0,T]$ to $L^\infty(\Omega;\Rt)\times L^p(\Omega;\Rtt)$ is measurable with respect to the class of Borel sets in the codomain, where $L^\infty(\Omega;\Rt)$ and $L^p(\Omega;\Rtt)$ are equipped with the weak-* topology and the weak topology, respectively, see, e.g. \cite[Section 6]{dalmaso.lazzaroni}. 

Also, similarly to \cite[Proposition 9.3]{mora.corral}, 
it is possible to show that there exists a solution $t\mapsto (C(t),K(t),\boldsymbol{y}(t))$ for which
\begin{equation*}
	C(t)=C_0 \cup \bigcup_{s\in [0,t]} C_{\boldsymbol{y}(s)} \quad \text{for all $t\in [0,T]$.}
\end{equation*}

  We will not present the proof for the last two statements as the  arguments  are identical to the ones provided in the references.

\section{Preliminaries}
\label{sec:preliminaries}

 In this section,  we collect  some results taken from the literature.

\subsection{Hausdorff distance}
  On the collection $\mathcal{K}(\closure{\Omega})$ of all compact subsets of $\closure{\Omega}$, \EEE  we consider the Hausdorff distance
\begin{equation*}
	d(K,\widetilde{K})\coloneqq \max\left\{ \sup_{\widetilde{\boldsymbol{x}}\in \widetilde{K}} \dist(\widetilde{\boldsymbol{x}};K), \sup_{\strut \boldsymbol{x}\in K} \dist(\boldsymbol{x};\widetilde{K})  \right\}
	 \quad \text{for all $K,\widetilde{K}\in\mathcal{K}(\closure{\Omega})$.}
\end{equation*} 
 which    induces  a complete metric on $\mathcal{K}(\closure{\Omega})$.   In particular,  $\mathcal{K}_m(\closure{\Omega})$ is a closed subset of $\mathcal{K}(\closure{\Omega})$ for every $m\in\N$, as shown in    \cite[Corollary 3.3]{dalmaso.toader}.   We write $K_n\to K$ in $\mathcal{K}(\closure{\Omega})$ whenever $d(K_n,{K}) \to 0$, as $n \to \infty$.  
The class $\mathcal{C}(\closure{\Omega})$  of finite subsets of $\closure{\Omega}$ clearly satisfies $\mathcal{C}(\closure{\Omega})\subset \mathcal{K}(\closure{\Omega})$. Hence, on $\mathcal{X}_m(\closure{\Omega})$,   we can consider the product metric
\begin{equation*}
	d\left( (C,K),(\widetilde{C},\widetilde{K})   \right)\coloneqq d(C,\widetilde{C})+d(K,\widetilde{K}) \quad \text{for all  $(C,K),(\widetilde{C},\widetilde{K}) \in \mathcal{X}_m(\closure{\Omega})$.}
\end{equation*}
 We write $ (C_{n},K_{n})\to (C,K)$ in $\mathcal{X}_m(\closure{\Omega})$ whenever $d( (C_n,K_n), (C,K)) \to 0$, as $n \to \infty$.  In general, the class $\mathcal{C}(\closure{\Omega})$ is not closed, however limits of sequences of sets in $\mathcal{C}(\closure{\Omega})$ with uniformly bounded cardinality belong to the same set. This follows from the following result concerning the Hausdorff convergence for finite sets which relates to the notion of componentwise convergence given  in   \cite[Definition 6.1]{mora.corral}.

\begin{lemma}[Hausdorff convergence for finite sets]\label{lem:hausdorff-finite}
Let $(C_n)_{n}$ be a sequence in $\mathcal{C}(\closure{\Omega})$ satisfying
\begin{equation*}
	\sup_ {n\in \N}\H^0(C_n)<+\infty.	
\end{equation*}
Let $C\in \mathcal{K}(\closure{\Omega})$ and suppose that $C_n\to C$ in $\mathcal{K}(\closure{\Omega})$.  Then, there exists a subsequence $(C_{n_k})_{k}$ and  two finite sets of indices $J\subset I \subset \N$ with $\H^0(C_{n_{k}})=\H^0(I)$ for all $k\in \N$ and  $\H^0(C)=\H^0(J)$ for which we can write $C_{n_{k}}=\{ \boldsymbol{a}^{n_{k}}_i:\:i\in I  \}$ and  $C=\{\boldsymbol{a}_j:\:j\in J\}$,   in such a way that:
 \begin{enumerate}[(i)]
 	\item For each $i\in I$, there exists a unique $j\in J$ for which $\boldsymbol{a}^{n_{k}}_i\to \boldsymbol{a}_j$, as $k\to\infty$;
 	\item For each $j\in J$, there exists $i\in I$ for which $\boldsymbol{a}^{n_{k}}_i\to \boldsymbol{a}_j$, as $k\to\infty$.
 \end{enumerate}  
Moreover, there holds
 \begin{equation}
 	\label{eq:lsc-hausdorff-finite}
 	\H^0(C)\leq \liminf_{n\to \infty} \H^0(C_n),
 \end{equation}
  so that $C\in\mathcal{C}(\closure{\Omega})$. 
\end{lemma}
\begin{proof}
As the convergence with respect to the Hausdorff distance is equivalent to the convergence in the sense of Kuratowski \cite[Proposition 4.4.14]{ambrosio.tilli}, we have the following:
\begin{enumerate}[(a)]
	\item For each sequence $(\boldsymbol{x}_n)_n\in \prod_n C_n$, for each  subsequence  $(\boldsymbol{x}_{n_k})_k$, and for each $\boldsymbol{x}\in \closure{\Omega}$ such that $\boldsymbol{x}_{n_k}\to \boldsymbol{x}$, as $k\to \infty$, there holds $\boldsymbol{x}\in C$;
	\item For each  $\boldsymbol{x}\in C$, there exists $(\boldsymbol{x}_n)_n\in \prod_n C_n$ such that $\boldsymbol{x}_n\to \boldsymbol{x}$, as $n\to \infty$.
\end{enumerate}

Here and below, $\prod_n C_n$ denotes the  collection of   sequences $(\boldsymbol{x}_n)_n$ in $\closure{\Omega}$ satisfying $\boldsymbol{x}_n\in C_n$ for all $n\in \N$.  

 By assumption,  the sequence $ \big  (\H^0(C_{n}) \big )_n$ is bounded and, hence, admits a subsequence converging to some $\beta\in \N$. Without relabeling, we assume that  $\H^0(C_n)=\beta$ for every $n\in \N$. Set $I\coloneqq \{1,\dots,\beta  \}$. \EEE   
We write $C_{n}=\{\boldsymbol{a}^{n}_i:\,i\in I  \}$ for all $n\in \N$. Since  $\closure{\Omega}$ is compact, we find a subsequence indexed by $(n_{k})_k$ and a finite set of points $\widehat{C}\coloneqq\{{\boldsymbol{a}}_i:\,i\in I\}\subset  \closure{\Omega}$   such that $\boldsymbol{a}^{n_{k}}_i\to {\boldsymbol{a}}_i$, as $k\to \infty$, for all $i\in I$.
By (a), we have $\widehat{C}\subset  C$. By removing repeated elements, we can write $\widehat{C}=\{\boldsymbol{a}_j:\,j\in J  \}$ for some $J\subset I$ with $\H^0(\widehat{C})=\H^0(J)$.  This  proves (i).  

We claim that  $C=\widehat{C}$.  
By contradiction, let ${\boldsymbol{x}}\in C\setminus \widehat{C}$ and  $0<r
 <  \frac{1}{2}  \min \{|{\boldsymbol{x}}-{\boldsymbol{a}}_j|: \:j\in J \}$.
By (b), there exists $({\boldsymbol{x}}_{n})_n\in \prod_n C_{n}$ such that ${\boldsymbol{x}}_{n}\to {\boldsymbol{x}}$, as $n\to \infty$. In particular, ${\boldsymbol{x}}_{n}\in B({\boldsymbol{x}},r)$ for $n\gg 1$. However, as ${\boldsymbol{x}}_{n}\in C_{n}$ for all $n\in \N$,  we also   have ${\boldsymbol{x}}_{n}\in \bigcup_{j\in J} B({\boldsymbol{a}}_j,r)$ for $n\gg 1$. Given our choice of $r$, this provides a contradiction. Therefore, $C=\widehat{C}$.

Let $j\in J$. Again by (b), there exists $({\boldsymbol{x}}_{n})_n\in \prod_n C_{n}$ such that ${\boldsymbol{x}}_{n}\to {\boldsymbol{a}}_j$, as $n\to \infty$. Hence, we find a sequence $(i_n)_n$ in $I$ such that $\boldsymbol{x}_n=\boldsymbol{a}^n_{i_n}$ for all $n\in \N$. Up to extracting a further subsequence  $(n_k)_k$,  we can assume that $i_{n_k}=i$ for some $i\in I$ and for all $k\in \N$. Thus, also (ii) is proved. 
 
Eventually, we check  \eqref{eq:lsc-hausdorff-finite}. 
Since the sequence $(n_k)_k$ can be chosen so that  $( \H^0(C_{n_k}) )_k$     converges to the right-hand side of \eqref{eq:lsc-hausdorff-finite}, we obtain
 \begin{equation*}
 	\H^0(C)=\H^0(J)\leq \H^0(I)=\lim_{k\to\infty} \H^0(C_{n_{k}})=\liminf_{n\to \infty} \H^0(C_n).
 \end{equation*}
\end{proof}

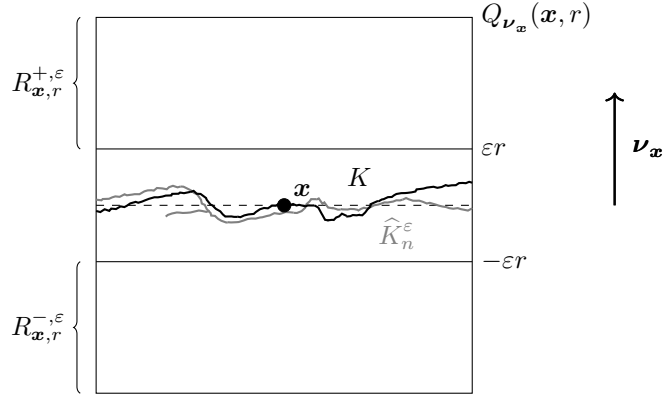
\begin{figure}
	\begin{tikzpicture}[scale=2.5]
		\draw (-1,-1) rectangle (1,1);
		\node[right] at (1,1) {$Q_{\boldsymbol{\nu}_{\boldsymbol{x}}}(\boldsymbol{x},r)$};
		\draw[dashed] (-1,0)--(1,0);
		\node[above right] at (0,0) {$\boldsymbol{x}$};
		\draw (-1,.3) -- (1,.3); \draw (-1,-.3)--(1,-.3);
		\draw[decoration={brace,raise=6pt},decorate] (-1,.3) -- node[left=8pt] {$R^{+,\varepsilon}_{\boldsymbol{x},r}$} (-1,1);
		\draw[decoration={brace,raise=6pt,mirror},decorate] (-1,-.3) -- node[left=8pt] {$R^{-,\varepsilon}_{\boldsymbol{x},r}$} (-1,-1);
		\draw[xshift=50pt,->,thick] (0,0) -- node[right=3pt] {$\boldsymbol{\nu}_{\boldsymbol{x}}$} (0,.6);
		\node[right] at (1,.3) {$\varepsilon r$};
		\draw[xshift=50pt,->,thick] (0,0) -- node[right=3pt] {$\boldsymbol{\nu}_{\boldsymbol{x}}$} (0,.6);
		\node[right] at (1,-.3) {$-\varepsilon r$};
		\draw [color=gray, thick, decorate, decoration={random steps,segment length=1.5pt,amplitude=.6pt}]
		(-1,.03)   to[out=0,in=180,distance=1mm] (-.55,.1);
		\draw [color=gray, thick, decorate, decoration={random steps,segment length=1.5pt,amplitude=.6pt}]
		(-.55,.1)   to[out=0,in=190,distance=1mm] (-.3,-.09);
		\draw [color=gray, thick, decorate, decoration={random steps,segment length=1.5pt,amplitude=.4pt}]
		(-.3,-.09)   to[out=0,in=170,distance=.8mm] (.05,-0.04);
		\draw [color=gray, thick, decorate, decoration={random steps,segment length=1.5pt,amplitude=.4pt}]
		(.05,-0.04)   to[out=0,in=180,distance=.7mm] (.17,.03);
		\draw [color=gray, thick, decorate, decoration={random steps,segment length=1.3pt,amplitude=.8pt}]
		(.17,.03)   to[out=350,in=200,distance=.4mm] (.3,-.02);
		\draw [color=gray, thick, decorate, decoration={random steps,segment length=1.3pt,amplitude=.6pt}]
		(.3,-.02)   to[out=0,in=130,distance=.4mm] (.39,-.03);
		\draw [color=gray, thick, decorate, decoration={random steps,segment length=1.3pt,amplitude=.4pt}]
		(.39,-.03)   to[out=30,in=180,distance=.2mm] (.63,.03);
		\draw [color=gray, thick, decorate, decoration={random steps,segment length=1.3pt,amplitude=.4pt}]
		(.63,.03)   to[out=10,in=200,distance=.55mm] (1,-.02);
		\draw [color=gray, thick, decorate, decoration={random steps,segment length=1.4pt,amplitude=.4pt}]
		(-.63,-.06)   to[out=30,in=0,distance=.7mm] (-0.4,-.03);
		\node[color=gray] at (.6,-.15) {$\widehat{K}^\varepsilon_n$};
		\draw[fill=black] (0,0) circle (1pt);
		\node at (.4,.15) {${K}$};
		\draw [ thick, decorate, decoration={random steps,segment length=1.5pt,amplitude=.6pt}]
		(-1,-.03)   to[out=0,in=180,distance=1mm] (-.5,.07);
		\draw [ thick, decorate, decoration={random steps,segment length=1.5pt,amplitude=.6pt}]
		(-.5,.07)   to[out=0,in=190,distance=1mm] (-.25,-.06);
		\draw [ thick, decorate, decoration={random steps,segment length=1.5pt,amplitude=.4pt}]
		(-.25,-.06)   to[out=0,in=170,distance=.8mm] (0,0);
		\draw [ thick, decorate, decoration={random steps,segment length=1.5pt,amplitude=.6pt}]
		(0,0)   to[out=0,in=0,distance=.4mm] (.16,0);
		\draw [ thick, decorate, decoration={random steps,segment length=1.3pt,amplitude=.8pt}]
		(.16,0)   to[out=350,in=200,distance=.4mm] (.24,-.07);
		\draw [ thick, decorate, decoration={random steps,segment length=1.3pt,amplitude=.6pt}]
		(.23,-.08)   to[out=0,in=130,distance=.4mm] (.33,-.062);
		\draw [ thick, decorate, decoration={random steps,segment length=1.3pt,amplitude=.4pt}]
		(.33,-.06)   to[out=10,in=90,distance=.2mm] (.4,-.06);
		\draw [ thick, decorate, decoration={random steps,segment length=1.3pt,amplitude=.4pt}]
		(.4,-.052)   to[out=10,in=180,distance=.55mm] (.5,.02);
		\draw [ thick, decorate, decoration={random steps,segment length=1.3pt,amplitude=.4pt}]
		(.5,.02)   to[out=40,in=180,distance=.7mm] (1,.12);
	\end{tikzpicture}
	\caption{  Representation of the statement of Lemma \ref{lem:separation}. The portions of both cracks $\widehat{K}^\varepsilon_n$ and $K$ inside of $Q_{\boldsymbol{\nu}_{\boldsymbol{x}}}(\boldsymbol{x},r)$ are contained in the strip $Q_{\boldsymbol{\nu}_{\boldsymbol{x}}}(\boldsymbol{x},r)\setminus (R^{+,\varepsilon}_{\boldsymbol{x},r}\cup R^{-,\varepsilon}_{\boldsymbol{x},r} )$. }
	\label{fig:R}
\end{figure}

The next result has been proved in \cite[ Theorem 3.2(c),  Proposition 3.3, and Remark 3.4]{francfort.giacomini.lopez}. 
\begin{lemma}[Local separation property of cracks] 
	\label{lem:separation}
Let $(K_n)_n$ be a sequence in $\mathcal{K}_1(\closure{\Omega})$ and let $K\in \mathcal{K}_1(\closure{\Omega})$ be such that $K_n \to K$ in $\mathcal{K}(\closure{\Omega})$.		Then, there exists a set $N\subset K$ with $\H^1(N)=0$  such that the following holds:   for all $\boldsymbol{x}\in K\setminus N$ and for all $\varepsilon>0$, there exist $\boldsymbol{\nu}_{\boldsymbol{x}}\in \So$, $ \bar{r}_{\boldsymbol{x}}^\varepsilon\EEE >0$, and $\bar{n}_{\boldsymbol{x}}^\varepsilon\EEE \in\N$ such that, for each $r\in \big(0,\bar{r}_{\boldsymbol{x}}^\varepsilon \big)$, there exists a sequence $(\widehat{K}_n^\varepsilon)_n$ in $ \mathcal{K}_1(\closure{\Omega})$ satisfying
\begin{equation*}
	K_n\subset \widehat{K}_n^\varepsilon,\quad \widehat{K}_n^\varepsilon \setminus K_n \subset Q_{\boldsymbol{\nu}_{\boldsymbol{x}}}(\boldsymbol{x},r),\quad \haus(\widehat{K}_n^\varepsilon \setminus K_n)\leq 3 \varepsilon r \qquad \text{for all $n\in \N$}
\end{equation*}
such that the two sets
\begin{align}
	\label{eq:R+}
	R^{+,\varepsilon}_{\boldsymbol{x},r}&\coloneqq \left\{ \boldsymbol{z}\in Q_{\boldsymbol{\nu}_{\boldsymbol{x}}}(\boldsymbol{x},r):\: (\boldsymbol{z}-\boldsymbol{x})\cdot \boldsymbol{\nu}_{\boldsymbol{x}}>\varepsilon r \right\},\\ 
	\label{eq:R-}
	R^{-,\varepsilon}_{\boldsymbol{x},r}&\coloneqq \left\{ \boldsymbol{z}\in Q_{\boldsymbol{\nu}_{\boldsymbol{x}}}(\boldsymbol{x},r):\: (\boldsymbol{z}-\boldsymbol{x})\cdot \boldsymbol{\nu}_{\boldsymbol{x}}<-\varepsilon r \right\}
\end{align}
belong to different connected components of $Q_{\boldsymbol{\nu_{\boldsymbol{x}}}}(\boldsymbol{x},r)\setminus K$  and  $Q_{\boldsymbol{\nu_{\boldsymbol{x}}}}(\boldsymbol{x},r)\setminus \widehat{K}^\varepsilon_n$ for all $n\geq  \bar{n}_{\boldsymbol{x}}^\varepsilon$ (see Figure \ref{fig:R}).
\end{lemma}

In the previous lemma, the sequence   $(\widehat{K}^\varepsilon_n)_n$  depends on both $\boldsymbol{x}$ and $r$ but, for simplicity, we omit this dependence  in  the notation. 

\subsection{Modeling of cavitation}
\label{subsec:cavitation}
In this subsection, we recall some facts about the variational theory of cavitation in nonlinear elasticity \cite{henao.moracorral.lusin}. 

  We   start   by introducing the topological degree,  referring to  \cite{fonseca.gangbo}  for a comprehensive account on this notion. Given $K\in\mathcal{K}(\closure{\Omega})$ and $\boldsymbol{y}\in W^{1,p}(\Omega\setminus K;\Rt)$, its precise representative $\boldsymbol{y}^*\colon \Omega \setminus K\to \Rt$ is defined as
\begin{equation*}
	\boldsymbol{y}^*(\boldsymbol{x})\coloneqq \limsup_{r\to 0^+} \dashint_{B(\boldsymbol{x},r)} \boldsymbol{y}(\boldsymbol{z})\,\d\boldsymbol{z},
\end{equation*}
while the set $L_{\boldsymbol{y}}$ of Lebesgue points of $\boldsymbol{y}$ is formed by those $\boldsymbol{x}\in \Omega \setminus K$ for which
\begin{equation*}
	\lim_{r \to 0^+} \dashint_{B(\boldsymbol{x},r)} |\boldsymbol{y}(\boldsymbol{z})-\boldsymbol{y}^*(\boldsymbol{x})|\,\d\boldsymbol{z}=0.
\end{equation*}
It is known that $\H^\alpha(\Omega \setminus (K \cup L_{\boldsymbol{y}}))=0$ for all $\alpha> (2-p)$. In particular, given  $\boldsymbol{a}\in \Omega \setminus K$, for  almost all   $r\in ( 0,\dist(\boldsymbol{a};\partial \Omega \cup K) )$ there holds $S(\boldsymbol{a},r)\subset L_{\boldsymbol{y}}$. 

\begin{definition}[Topological degree and topological image]
Let $K\in\mathcal{K}(\closure{\Omega})$ and $\boldsymbol{y}\in W^{1,p}(\Omega\setminus K;\Rt)$. Given $\boldsymbol{a}\in \Omega \setminus K$ and $r\in ( 0,\dist(\boldsymbol{a};\partial \Omega \cup K))$ such that $S(\boldsymbol{a},r)\subset L_{\boldsymbol{y}}$ and $\boldsymbol{y}\in W^{1,p}(S(\boldsymbol{a},r);\Rt)$, we define the \emph{topological degree} of $\boldsymbol{y}$ with respect to $B(\boldsymbol{a},r)$ to be the topological degree of any extension of $\boldsymbol{y}^*\restr{S(\boldsymbol{a},r)}$ with respect to  $B(\boldsymbol{a},r)$  and we denote it as $\deg(\boldsymbol{y},B(\boldsymbol{a},r),\cdot)\colon \Rt\setminus \boldsymbol{y}^*(S(\boldsymbol{a},r)) \to \mathbb{Z}$. Then, we define the \emph{topological image} of $B(\boldsymbol{a},r)$ via $\boldsymbol{y}$ as
\begin{equation*}
	\imt(\boldsymbol{y},B(\boldsymbol{a},r))\coloneqq \left\{ \boldsymbol{\xi}\in\Rt \setminus \boldsymbol{y}^*(S(\boldsymbol{a},r)):\:\deg(\boldsymbol{y},B(\boldsymbol{a},r),\boldsymbol{\xi})\neq 0  \right\}.
\end{equation*} 
\end{definition} 
Note that $S(\boldsymbol{a},r)\subset L_{\boldsymbol{y}}$ and $\boldsymbol{y}\in W^{1,p}(S(\boldsymbol{a},r);\Rt)$ imply $\boldsymbol{y}^*\in C^0(S(\boldsymbol{a},r);\Rt)$ by Morrey's embedding. This together with   the exclusive dependence of the degree on the values at the boundary   ensures  that $\deg(\boldsymbol{y},B(\boldsymbol{a},r),\cdot)$ is well defined. By the continuity of the topological degree, the set $\imt(\boldsymbol{y},B(\boldsymbol{a},r))$ is open and bounded with  $\partial\,\imt(\boldsymbol{y},B(\boldsymbol{a},r))\subset\boldsymbol{y}^*(S(\boldsymbol{a},r))$. 

Next, we provide a rigorous definition for cavities and cavitation points.

\begin{definition}[Topological image of a point and cavitation points]
	\label{def:cavitation-points}
Let $K\in\mathcal{K}(\closure{\Omega})$ and $\boldsymbol{y}\in W^{1,p}(\Omega\setminus K;\Rt)$. Given $\boldsymbol{a}\in \Omega \setminus K$, we define the \emph{topological image of the point} $\boldsymbol{a}$ via $\boldsymbol{y}$ as
\begin{equation*}
	\imt(\boldsymbol{y},\boldsymbol{a})\coloneqq \bigcap  \left\{ \closure{\imt(\boldsymbol{y},B(\boldsymbol{a},r))}:\: \text{$r\in \big( 0,\dist(\boldsymbol{a};\partial \Omega \cup K) \big)$ s.t. $S(\boldsymbol{a},r)\subset L_{\boldsymbol{y}}$ and $\boldsymbol{y}\in W^{1,p}(S(\boldsymbol{a},r);\Rt)$}  \right\}.
\end{equation*}	
Moreover, the set of \emph{cavitation points} is defined as
\begin{equation*}
	C_{\boldsymbol{y}}\coloneqq \left \{ \boldsymbol{a}\in \Omega \setminus K:\:\leb(\imt(\boldsymbol{y},\boldsymbol{a}))>0   \right \}.
\end{equation*}
\end{definition}
As observed in \cite[Lemma 7.3(i)]{mueller.spector}, one has $\closure{\imt(\boldsymbol{y},B(\boldsymbol{a},r_1))}\subset \closure{\imt(\boldsymbol{y},B(\boldsymbol{a},r_2))}$ for $r_1,r_2 \in \big( 0,\dist(\boldsymbol{a};\partial \Omega \cup K) \big)$ as above with $r_1<r_2$. Thus, the set $\imt(\boldsymbol{y},\boldsymbol{a})$ is given by a decreasing intersection of  nonempty  compact sets. As such, $\imt(\boldsymbol{y},\boldsymbol{a})$ is  also   compact and nonempty. 

Next, we recall  the invertibility condition introduced in  \cite{mueller.spector}. In this paper, we will often refer to a localized version of this condition as given below. \EEE 

\begin{definition}[Invertibility condition]
	\label{def:INV}
  Let $K\in\mathcal{K}(\closure{\Omega})$ and $\boldsymbol{y}\in W^{1,p}(\Omega\setminus K;\Rt)$. Given an open set $U\subset \Omega \setminus K$,  the map $\boldsymbol{y}$ is termed to satisfy condition (INV)   in $U$   whenever  for all $\boldsymbol{a}\in U$ and for almost all $r \in \big (0,\dist(\boldsymbol{a};\partial U)\big) $ the following  holds: 
\begin{enumerate}[(i)]
	\item $\boldsymbol{y}(\boldsymbol{x})\in \imt(\boldsymbol{y}, B(\boldsymbol{a},r)  )$ for almost every $\boldsymbol{x}\in B(\boldsymbol{a},r)$;
	\item $\boldsymbol{y}(\boldsymbol{x})\notin \imt(\boldsymbol{y}, B(\boldsymbol{a},r)  )$ for almost every $\boldsymbol{x}\in U \setminus \closure{B}(\boldsymbol{a},r)$.
\end{enumerate}
\end{definition}

By \cite[Lemma 3.4]{mueller.spector}, if $\det D\boldsymbol{y}\neq 0$ almost everywhere in $U$ and $\boldsymbol{y}$ satisfies (INV) in $U$, then $\boldsymbol{y}$ is almost everywhere injective in $U$ meaning that $\boldsymbol{y}\restr{U \setminus N}$ is injective for some set $N\subset U$ with $\leb(N)=0$.
Note that $U=\Omega\setminus K$ is an admissible choice for Definition \ref{def:INV}.   

Henceforth, \emph{given $K\in \mathcal{K}(\closure{\Omega})$ and $\boldsymbol{y}\in W^{1,p}(\Omega\setminus K;\Rt)$, we systematically identify $\boldsymbol{y}$, $D\boldsymbol{y}$, and $\det D \boldsymbol{y}$ with their extensions to $\Omega$ by zero}.  The next result from \cite{mueller.spector} concerns the stability of condition (INV)  for  the   convergence of cracks with respect to the Hausdorff distance and the weak convergence of deformations.

\begin{lemma}[Stability of the invertibility condition]
	\label{lem:INV}
Let $(K_n)_{ n}$ be a sequence in $\mathcal{K}(\closure{\Omega})$ and $K\in \mathcal{K}(\closure{\Omega})$. Also, let $(\boldsymbol{y}_n)_n$ be a sequence with $\boldsymbol{y}_n\in W^{1,p}(\Omega \setminus K_n;\Rt)$ satisfying condition {\rm (INV)}  in $\Omega\setminus K_n$  for all $n\in \N$ and $\boldsymbol{y}\in W^{1,p}(\Omega \setminus K ;\Rt)$. Suppose that 
\begin{equation*}
	\text{$K_n\to K$ in $\mathcal{K}(\closure{\Omega})$}, \qquad \text{$\boldsymbol{y}_n \wk \boldsymbol{y}$ in $L^p(\Omega;\Rt)$}, \qquad \text{$D\boldsymbol{y}_n \wk D\boldsymbol{y}$ in $L^p(\Omega;\Rtt)$.}
\end{equation*}
Then, $\boldsymbol{y}$ satisfies condition {\rm (INV)}  in $\Omega \setminus K$. 
\end{lemma}
\begin{proof}
Let $\boldsymbol{a}\in \Omega\setminus K$ and $U\subset \subset \Omega \setminus K$ be open with $\boldsymbol{a}\in U$. By Hausdorff convergence, $U\subset \Omega \setminus K_n $ for $n\gg 1$ depending on $U$. By applying  \cite[Lemma 3.3]{mueller.spector}   to $(\boldsymbol{y}_n\restr{U})_{n\gg 1}$ and $\boldsymbol{y}\restr{U}$, we see that   $\boldsymbol{y}$ satisfies  condition (INV) in $U$.  Taking a sequence $(U_j)_j$ with $U_j\subset \subset \Omega \setminus K$ whose union coincides with  $\Omega \setminus K$,  we obtain that $\boldsymbol{y}$ satisfies condition (INV) on the whole set $\Omega\setminus K$.
\end{proof}

We  recall   the surface energy functional introduced in \cite{henao.moracorral.invertibility}.   We refer to \cite{henao.moracorral.fracture} for the characterization of this energy in geometric terms. Here, we present it in a localized version.  

\begin{definition}[Surface energy]
	\label{def:surface-energy}
Let $K\in \mathcal{K}(\closure{\Omega})$ and  $\boldsymbol{y}\in W^{1,p}(\Omega\setminus K;\Rt)$.  Given an open set $U\subset \Omega \setminus K$, we define the   corresponding   surface energy functional
\begin{equation*}
	\mathcal{S}_U(\boldsymbol{y})\coloneqq \sup \left\{ \mathcal{S}^{\boldsymbol{y}}(\boldsymbol{\eta}):\: \boldsymbol{\eta}\in C^\infty_{\rm c}((\Omega \setminus K) \times \Rt;\Rt),\: \supp \,\boldsymbol{\eta}\subset U  \times \R^2,  \:\: \|\boldsymbol{\eta}\|_{L^\infty(U\times \Rt;\Rt)}\leq 1     \right\},
\end{equation*}
where 
\begin{equation*}
	\mathcal{S}^{\boldsymbol{y}}(\boldsymbol{\eta})\coloneqq \int_{\Omega \setminus K} \left( \cof D \boldsymbol{y}(\boldsymbol{x}) : D_{\boldsymbol{x}}\boldsymbol{\eta}(\boldsymbol{x},\boldsymbol{y}(\boldsymbol{x}))+\div_{\boldsymbol{\xi}}\boldsymbol{\eta}(\boldsymbol{x},\boldsymbol{y}(\boldsymbol{x}))\det D \boldsymbol{y}(\boldsymbol{x})  \right)\,\d\boldsymbol{x}
\end{equation*}
for every function $\boldsymbol{\eta}\in C^\infty_{\rm c}((\Omega \setminus K) \times \Rt;\Rt)$ with variables $(\boldsymbol{x},\boldsymbol{\xi})\in (\Omega \setminus K) \times \Rt$.  
\end{definition}

Also here, the choice $U=\Omega \setminus K$ is admissible for the previous definition. In this case, we simply write $\mathcal{S}(\boldsymbol{y})\coloneqq \mathcal{S}_{\Omega \setminus K}(\boldsymbol{y})$, which coincides with the definition in \cite{henao.moracorral.invertibility} in the case of deformations defined in $\Omega\setminus K$.

  It has been proved in \cite{henao.moracorral.lusin} that, for  deformations  having finite surface energy, this energy equals the sum of the perimeters of the cavities created by the deformation. In our context, the result reads as follows.   

\begin{lemma}[Representation of the surface energy]
	\label{lem:EU}
Let $K\in \mathcal{K}(\closure{\Omega})$ and $\boldsymbol{y}\in W^{1,p}(\Omega \setminus K;\Rt)$ with $\mathcal{S}(\boldsymbol{y})<+\infty$  satisfying   {\rm (INV)} in $\Omega\setminus K$.  Then, for every open set $U\subset \Omega \setminus K$, we have 
\begin{equation*}
	\mathcal{S}_U(\boldsymbol{y})=\sum_{\boldsymbol{a}\in C_{\boldsymbol{y}}\cap U}
	 \per \left( \imt(\boldsymbol{y},\boldsymbol{a})  \right).
\end{equation*}
In particular,  recalling \eqref{eq:WPL}, we have   $\mathcal{P}(C,K,\boldsymbol{y})=\mathcal{S}(\boldsymbol{y})$ for every $(C,K,\boldsymbol{y})\in \mathcal{A}_{p,M}^m$.
\end{lemma}
 
As an immediate consequence of this result, we deduce that $C_{\boldsymbol{y}}$ is at most countable whenever $\mathcal{S}(\boldsymbol{y})<+\infty$.

\begin{proof}
It is sufficient to apply \cite[Theorem 4.6]{henao.moracorral.lusin} to $\boldsymbol{y}\restr{U}$. Note that the result still holds true for deformations in $W^{1,p}(U;\Rt)$ not necessarily bounded, as already observed in \cite[Proposition 2.14]{henao.moracorral.regularity} and \cite[Theorem 2.39]{bresciani.friedrich.moracorral}.
\end{proof}

\section{Proof of the main result}
\label{sec:proof}
 
This section is devoted to the proof of Theorem \ref{thm:existence}.

\subsection{Coercivity and lower semicontinuity}
In this subsection, we discuss the main properties of the energy functional $\mathcal{F}$ in \eqref{eq:F}--\eqref{eq:WPL}. 

  First, we show that the class $\mathcal{A}_{p,M}^m$ of admissible   states  is closed with respect to the relevant notion of convergence.  

\begin{lemma}[Closedness of the state space]
	\label{lem:closedness-state-space}
Let $  ((C_n,K_n,\boldsymbol{y}_n)  )_n$ be a sequence in $\mathcal{A}_{p,M}^m$ satisfying 
\begin{equation*}
		\label{eq:l0}
		\sup_{n\in \N}  \left\{ \mathcal{P}(C_n,K_n,\boldsymbol{y}_n) + \H^0(C_n)  \right\}<+\infty.
\end{equation*}
 Suppose that there exist    $C,K\in \mathcal{K}(\closure{\Omega})$,  $\boldsymbol{y}\in L^p(\Omega;\Rt)$, $\boldsymbol{F}\in L^p(\Omega;\Rtt)$, and $h\in L^1_+(\Omega)$   such that 
\begin{equation}
		\label{eq:l1}
		\text{$C_n\to C$ in $\mathcal{K}(\closure{\Omega})$,} \qquad \text{$K_n\to K$ in $\mathcal{K}(\closure{\Omega})$,} \qquad \text{$\boldsymbol{y}_n \wks \boldsymbol{y}$ in $L^\infty(\Omega;\Rt)$,}
\end{equation}
and
\begin{equation}
		\label{eq:l2}
		\text{$D\boldsymbol{y}_n \wk \boldsymbol{F}$ in $L^p(\Omega;\Rtt)$,} \qquad \text{$\det D \boldsymbol{y}_n\wk h$ in $L^1(\Omega)$.}
\end{equation}
Then,   $(C,K,\boldsymbol{y})\in \mathcal{A}_{p,M}^m$ and there hold $\boldsymbol{F}=D\boldsymbol{y}$ and $h=\det D \boldsymbol{y}$ in $\Omega\setminus K$. 
\end{lemma}
\begin{proof}
   First, $(C,K)\in\mathcal{X}_m(\closure{\Omega})$ by Lemma \ref{lem:hausdorff-finite} and the closedness of $\mathcal{K}_m(\closure{\Omega})$. 	
Clearly, $\|\boldsymbol{y}\|_{L^\infty(\Omega;\Rt)}\leq M$. 	
Let $U\subset \subset \Omega\setminus K$ be open. As $K_n \to K$ in $\mathcal{K}(\closure{\Omega})$,  we have $U\subset \subset \Omega \setminus K_n$ for $n\gg 1$ depending on $U$. This inclusion together with \eqref{eq:l1}--\eqref{eq:l2} yields $\boldsymbol{y}\in W^{1,p}(U;\Rt)$ with $D\boldsymbol{y}=\boldsymbol{F}$ in $U$.   As $U$ is arbitrary, we actually have  $\boldsymbol{y}\in W^{1,p}(\Omega \setminus K;\Rt)$ with $D\boldsymbol{y}=\boldsymbol{F}$ in $\Omega \setminus K$. Also, $\boldsymbol{y}$ satisfies (INV) in $\Omega \setminus K$ by Lemma \ref{lem:INV}.  

 In view of Lemma \ref{lem:EU}, we have
\begin{equation*}
	\sup_{n\gg 1}\mathcal{S}_U(\boldsymbol{y}_n)=\sup_{n\gg 1} \sum_{\boldsymbol{a}\in C_{\boldsymbol{y}_n}\cap U} \per\left(\imt( \boldsymbol{y}_n ,\boldsymbol{a})  \right)  \leq \sup_{n\in \N} \mathcal{P}(  C_n,  K_n,\boldsymbol{y}_n)<+\infty.
\end{equation*} 
Note that \eqref{eq:l1}--\eqref{eq:l2} and the identity $\boldsymbol{F}=D\boldsymbol{y}$ in $\Omega\setminus K$ yield $\boldsymbol{y}_n\wk \boldsymbol{y}$ in $W^{1,p}(U;\Rt)$. Thus, $\cof D\boldsymbol{y}_n \wk \cof D\boldsymbol{y}$ in $L^p(U;\Rtt)$ and, by Rellich's theorem,  we have   $\boldsymbol{y}_n\to  \boldsymbol{y}$ almost everywhere in $U$ along a not relabeled subsequence. Then,   
by applying \cite[Theorem~3]{henao.moracorral.invertibility} to $(\boldsymbol{y}_n\restr{U})_n$ and $\boldsymbol{y}\restr{U}$, we deduce $h=\det D\boldsymbol{y}$ in $U$ and the inequality
\begin{equation*}
	\mathcal{S}_U(\boldsymbol{y})\leq \liminf_{n\to \infty} \mathcal{S}_U(\boldsymbol{y}_n)\leq \sup_{n\in \N} \mathcal{S}(\boldsymbol{y}_n)=\sup_{n\in \N} \mathcal{P}(  C_n,  K_n,\boldsymbol{y}_n)<+\infty.
\end{equation*}
From this, taking the supremum among all $U\subset \subset \Omega \setminus K$, we see that  $h=\det D \boldsymbol{y}$ in $\Omega\setminus K$ and  $\mathcal{S}(\boldsymbol{y})<+\infty$.

Eventually, we need to prove that $C_{\boldsymbol{y}}\subset C \cap (\Omega \setminus K)$. 
First, by Lemma \ref{lem:hausdorff-finite}, there exist a subsequence $(C_{n_k})_k$ and two finite sets $J\subset I\subset \N$
such that we can write $C_{n_k}=\{ \boldsymbol{a}^{n_k}_i: \:i\in I \}$ for all $k\in \N$ and $C=\{\boldsymbol{a}_j:\:j\in J \}$ in such a way that   claim (i) of Lemma \ref{lem:hausdorff-finite} holds true.  Let $\boldsymbol{a}\in C_{\boldsymbol{y}} \subset \Omega \setminus K$. We show that $\boldsymbol{a}\in C$. Let $U\subset \subset \Omega \setminus K$  be open  with $\boldsymbol{a}\in U$  and consider $n\gg 1$ so that $U\subset \Omega\setminus K_n$. By applying \cite[Lemma 3.11]{bresciani.friedrich.moracorral} to $(\boldsymbol{y}_n\restr{U})_n$ and $\boldsymbol{y}\restr{U}$, we find a set $I_U$ and, for each $i\in I_U$, a sequence $(\boldsymbol{a}^{n_k}_i)_k$ of points in $U$ such that $\{ \boldsymbol{a}_i^{n_k}:\:i\in I_U \}\subset C_{\boldsymbol{y}_{n_k}}\cap U$ for all $k\in\N$, and we can write $C_{\boldsymbol{y}} \cap U=\{ \widetilde{\boldsymbol{a}}_i:\:i\in {I}_U \}$ in such a way that $\boldsymbol{a}^{n_k}_i\to \widetilde{\boldsymbol{a}}_i$, as $k\to \infty$, for all $i\in I_U$. Given that $C_{\boldsymbol{y}_{n_k}}\subset C_{n_k}$, up to identification, we may assume $I_U\subset I$. Let $i\in I_U$ be such that $\boldsymbol{a}=\widetilde{\boldsymbol{a}}_i$. Then, $\boldsymbol{a}^{n_k}_i\to \boldsymbol{a}$  as $k\to \infty$, so that claim (i) of Lemma \ref{lem:hausdorff-finite} yields $\boldsymbol{a}=\boldsymbol{a}_j$ for some $j\in J$. Therefore, $\boldsymbol{a}\in C$. \EEE 
\end{proof}
 
 Next, we investigate  coercivity and lower semicontinuity of $\mathcal{F}$. 
 
\begin{proposition}[Compactness]
	\label{prop:coercivity}
Let $  ((C_n,K_n,\boldsymbol{y}_n)   )_n$ be a sequence in $\mathcal{A}^m_{p,M}$  satisfying
\begin{equation*}
	\sup_{n\in \N} \mathcal{F}(t,C_n,K_n,\boldsymbol{y}_n)<+\infty	
\end{equation*}
for some $t\in [0,T]$. Then, there exist   a not relabeled subsequence and \EEE  $(C,K,\boldsymbol{y})\in \mathcal{A}^m_{p,M}$ such that 
\begin{equation}
	\label{eq:l3}
	  \text{$(C_{n},K_{n})\to (C,K)$ in $\mathcal{X}_m(\closure{\Omega})$,} \qquad \text{$\boldsymbol{y}_{n} \wks \boldsymbol{y}$ in $L^\infty(\Omega;\Rt)$,}
\end{equation}
and
\begin{equation}
	\label{eq:l4}
	  \text{$D\boldsymbol{y}_{n} \wk D\boldsymbol{y}$ in $L^p(\Omega;\Rtt)$,} \qquad \text{$\det D \boldsymbol{y}_{n}\wk \Det D \boldsymbol{y}$ in $L^1(\Omega)$,}
\end{equation}
as  $n\to \infty$.  
\end{proposition}
\begin{proof}
Thanks to assumption (W2),  we have   
\begin{equation*}
	\sup_{n\in \N} \left\{ \|D \boldsymbol{y}_n\|_{L^p(\Omega;\Rtt)} + \|\gamma(\det D \boldsymbol{y}_n)\|_{L^1(\Omega)} + \mathcal{P}(C_n,K_n,\boldsymbol{y}_n)+\H^0(C_n)  +\H^1(K_n)    \right\}<+\infty.
\end{equation*}
From this bound, the confinement condition, the De la Vall\'{e} Poussin criterion, and Blaschke's selection theorem \cite[Theorem 4.4.15]{ambrosio.tilli}, we conclude that there exist a subsequence  together   with  $\boldsymbol{y}\in L^\infty(\Omega;\Rt)$, $\boldsymbol{F}\in L^p(\Omega;\Rtt)$, $h\in  L^1(\Omega)   $, $C\in \mathcal{C}(\closure{\Omega})$, and $K\in \mathcal{K}(\closure{\Omega})$ such that \eqref{eq:l1}--\eqref{eq:l2} hold true. By a standard contradiction argument based on \eqref{eq:gamma}, we see that $h>0$ almost everywhere.  Thus, by applying Lemma \ref{lem:closedness-state-space}, we deduce that $(C,K,\boldsymbol{y})\in \mathcal{A}_{p,M}^m$ and we identify $\boldsymbol{F}=D\boldsymbol{y}$ and $h=\det D \boldsymbol{y}$  in $\Omega \setminus K$.  
\end{proof}

\begin{proposition}[Lower  semicontinuity]
	\label{prop:lsc}
Let $  ( (C_n,K_n,\boldsymbol{y}_n)  )_n$  be a sequence in $\mathcal{A}_{p,M}^m$ and let $(C,K,\boldsymbol{y})\in \mathcal{A}^m_{p,M}$. Suppose that \eqref{eq:l3}--\eqref{eq:l4} hold. Then,
\begin{align}
	\label{eq:lsc1}
	\mathcal{W}(\boldsymbol{y})&\leq \liminf_{n\to \infty} \mathcal{W}(\boldsymbol{y}_n), \hspace{40pt} \mathcal{P}(C,K,\boldsymbol{y})\leq \liminf_{n\to \infty} \mathcal{P}(C_n,K_n,\boldsymbol{y}_n),\\
	\label{eq:lsc1-bis} 	\H^0(C)&\leq \liminf_{n\to \infty} \H^0(C_n), \hspace{45pt} 	\H^1(K)\leq \liminf_{n\to \infty} \H^1(K_n)
\end{align}
and
\begin{equation}
\label{eq:lsc2}
\mathcal{L}(t,\boldsymbol{y})=\lim_{n\to \infty} \mathcal{L}(t,\boldsymbol{y}_n) \quad \text{for all $t\in [0,T]$.}
\end{equation}
Therefore,
\begin{equation*}
	\mathcal{F}(t,C,K,\boldsymbol{y})\leq \liminf_{n\to \infty} \mathcal{F}(t,C_n,K_n,\boldsymbol{y}_n) \quad \text{for all $t\in [0,T]$.}
\end{equation*} 	
\end{proposition}
\begin{proof}
Thanks to assumption (W3), from \eqref{eq:l3}--\eqref{eq:l4}, we have the first inequality in \eqref{eq:lsc1} and \eqref{eq:lsc2}.
Next, let $U\subset \subset \Omega \setminus K$. By \eqref{eq:l3}, we have $U\subset \Omega\setminus K_n$ for $n\gg 1$.
 Note that $\boldsymbol{y}_n\to \boldsymbol{y}$ almost everwhere in $U$, $\cof D \boldsymbol{y}_n \wk \cof D \boldsymbol{y}$ in $L^p(U;\Rtt)$, and $\det D \boldsymbol{y}_n \wk \det D \boldsymbol{y}$ in $L^1(U)$ by \eqref{eq:l3}--\eqref{eq:l4} and Rellich's embedding. Using Lemma \ref{lem:EU} and applying \cite[Theorem 3]{henao.moracorral.invertibility} to  $(\boldsymbol{y}_n\restr{U})_{n\gg1}$ and $\boldsymbol{y}\restr{U}$,   we find
\begin{equation*}
	\sum_{\boldsymbol{a}\in C_{\boldsymbol{y}}\cap U} \per \left( \imt(\boldsymbol{y},\boldsymbol{a}) \right)=\mathcal{S}_U(\boldsymbol{y})\leq \liminf_{n\to \infty} \mathcal{S}_U(\boldsymbol{y}_n)\leq \liminf_{n\to \infty} \mathcal{S}(\boldsymbol{y}_n)=\liminf_{n\to \infty} \mathcal{P}(C_n,K_n,\boldsymbol{y}_n).
\end{equation*} 
Then, taking the supremum among all open sets $U\subset \subset \Omega \setminus K$, we deduce the second inequality in \eqref{eq:lsc1}.
By Lemma \ref{lem:hausdorff-finite}, we obtain the first inequality in \eqref{eq:lsc1-bis}. Eventually,  thanks to  Go\l{}\k{a}b's  theorem \cite[Corollary 3.3]{dalmaso.toader},  we obtain   the second part of  \eqref{eq:lsc1-bis}.  
\end{proof}

\subsection{Incremental minimization problem and compactness}
In this subsection, we   introduce     an  incremental minimization problem 
 related to a  discretization of the time interval.
 
Henceforth, we  identify each ordered set $(t_i)_{i=0,\dots,\nu}$, where $\nu\in \N$ and $0=t_0<t_1<\dots<t_\nu=T$, with the partition $\{[t_{i-1},t_i]   \}_{i=1,\dots,\nu}$ of $[0,T]$. Let $(t_i)_{i=0,\dots,\nu}$  be a partition of $[0,T]$  and let $(C_0,K_0,\boldsymbol{y}_0)\in\mathcal{A}^m_{p,M}$. The \emph{incremental minimization problem} determined by $(t_i)_{i=0,\dots,\nu}$ with initial datum $(C_0,K_0,\boldsymbol{y}_0)$ aims at finding $  ((C_i,K_i,\boldsymbol{y}_i)  )_{i=1,\dots,\nu} \in (\mathcal{A}_{p,M}^m)^\nu$ satisfying the  unilateral stability condition
\begin{equation}
	\label{eq:imp}
	\mathcal{F}(t_i,C_i,K_i,\boldsymbol{y}_i)\leq  \mathcal{F}(t_i,\widetilde{C},\widetilde{K},\widetilde{\boldsymbol{y}}) \quad \text{for all $(\widetilde{C},\widetilde{K},\widetilde{\boldsymbol{y}})\in\mathcal{A}^m_{p,M}$ with $\widetilde{C}\supset C_{i-1}$ and $\widetilde{K}\supset K_{i-1}$}
\end{equation}
for all $i=1,\dots,\nu$.

The next result  ensures  the existence of   solutions to   the incremental minimization problem and shows that these  satisfy  a discrete energy-dissipation inequality.

\begin{proposition}[Incremental minimization problem]
	\label{prop:imp}
Let $(t_i)_{i=0,\dots,\nu}$ be a partition of $[0,T]$ and  let  $(C_0,K_0,\boldsymbol{y}_0)\in\mathcal{A}^m_{p,M}$     satisfy \eqref{eq:initial-datum}. Then:
\begin{enumerate}[(i)]
	\item There exists a solution $ ( (C_i,K_i,\boldsymbol{y}_i) )_{i=1,\dots,\nu} \in(\mathcal{A}^m_{p,M})^\nu$ of the incremental minimization problem determined by $(t_i)_{i=0,\dots,\nu}$ with initial datum $(C_0,K_0,\boldsymbol{y}_0)$;
	\item For all $i= 1,\dots,\nu$, there holds 
	\begin{equation*}
		\mathcal{F}(t_i,C_i,K_i,\boldsymbol{y}_i)\leq \mathcal{F}(0,C_0,K_0,\boldsymbol{y}_0)-\sum_{j=1}^i \int_{t_{j-1}}^{t_j} \int_{\Omega} \dt{\boldsymbol{f}}(s)\cdot \boldsymbol{y}_{j-1}\,\d\boldsymbol{x}\,\d s.
	\end{equation*}
\end{enumerate}
\end{proposition}
\begin{proof}
(i) Let $i\in \{1,\dots,\nu\}$ and suppose that  $(C_{i-1},K_{i-1},\boldsymbol{y}_{i-1})$  is known. With the aid of Proposition \ref{prop:coercivity} and Proposition \ref{prop:lsc}, the existence of $(C_i,K_i,\boldsymbol{y}_i)\in\mathcal{A}_{p,M}^m$ satisfying \eqref{eq:imp} is established by the direct method of the calculus of variations. In this regard, we observe that, if $ (( C^n,K^n)  )_n$ is a sequence in $\mathcal{X}_m(\closure{\Omega})$ with $C^n\supset C_{i-1}$ and $K^n \supset K_{i-1}$ for every $n\in\N$ and there exists $(C,K)\in\mathcal{X}_m(\closure{\Omega})$ such that $(C^n,K^n)\to (C,K)$ in $\mathcal{X}_m(\closure{\Omega})$, then necessarily $C\supset C_{i-1}$ and $K \supset K_{i-1}$. 

(ii) The inequality is obtained by testing \eqref{eq:imp} with $(\widetilde{C},\widetilde{K},\widetilde{\boldsymbol{y}})=(C_{i-1},K_{i-1},\boldsymbol{y}_{i-1})$
 and summing over $i=1,\dots,\nu$ after cancellations in a telescopic sum. 
\end{proof}

We now consider a  sequence  $ ( (t^n_i)_{i=0,\dots,\nu_n}   )_n$   of partitions with asymptotically vanishing size,   that is
\begin{equation}
	\label{eq:size-partition}
	\Delta_n\coloneqq \max_{i=1,\dots,\nu_n }  |t_i^n-t_{i-1}^n|\to 0 \quad \text{as $n\to \infty$.}
\end{equation}
Let $(C_0,K_0,\boldsymbol{y}_0)\in \mathcal{A}^m_{p,M}$ and, for every $n\in\N$, let $ ( (C_i^n,K_i^n,\boldsymbol{y}_i^n)  )_{i=1,\dots,\nu_n}\in (\mathcal{A}^m_{p,M})^{\nu_n}$ \EEE be a solution of the incremental minimization problem determined by $(t^n_i)_{i=0,\dots,\nu_n}$ with initial datum $(C_0,K_0,\boldsymbol{y}_0)$. We define the corresponding piecewise-constant interpolations  $t\mapsto (C^n(t),K^n(t),\boldsymbol{y}^n(t))$ as the function from $[0,T]$ to $\mathcal{A}^m_{p,M}$ given by
\begin{equation}
	\label{eq:pc1}
	C^n(t)\coloneqq C_i^n,\quad K^n(t)\coloneqq K_i^n,\quad \boldsymbol{y}^n(t)\coloneqq \boldsymbol{y}_i^n \qquad \text{for $t\in [t_{i}^n,t_{i+1}^n)$ and $i= 0,  \dots,\nu_n-1$},
\end{equation}
 where $(C^n_0, K^n_0, \boldsymbol{y}_0^n)  = (C_0, K_0, \boldsymbol{y}_0) $,  and
\begin{equation}
	\label{eq:pc2}
	C^n(T)\coloneqq C^n_{\nu_n}, \quad K^n(T)\coloneqq K_{\nu_n}^n, \quad  \boldsymbol{y}^n(T)\coloneqq \boldsymbol{y}_{\nu_n}^n.
\end{equation}
Also,  we define the function $\tau^n\colon [0,T] \to \{t^n_0,\dots,t^n_{\nu_n}\}$ by setting
\begin{equation}
	\label{eq:tau}
	\text{$\tau^n(t)\coloneqq t^n_i$ for all $t\in [t^n_i,t^n_{i+1})$ and $i=  0,  \dots,\nu_n-1$, \qquad  $\tau^n(T)\coloneqq t^n_{\nu_n}=T$.}
\end{equation} 
\begin{proposition}[Compactness  for solutions of the incremental minimization problem]
	\label{prop:compactness}
Let $  ( (t^n_i)_{i=0,\dots,\nu_n}    )_n$ be a sequence of partitions of $[0,T]$ satisfying   \eqref{eq:size-partition}. \EEE Let $(C_0,K_0,\boldsymbol{y}_0)\in \mathcal{A}^m_{p,M}$ satisfy \eqref{eq:initial-datum}. 
 For every $n\in \N$, let   $ ( (C_i^n,K_i^n,\boldsymbol{y}_i^n)  )_{i=1,\dots,\nu_n}\in (\mathcal{A}^m_{p,M})^{\nu_n}$ \EEE be a solution of the incremental minimization problem determined by $(t^n_i)_{i=0,\dots,\nu_n}$ with initial datum $(C_0,K_0,\boldsymbol{y}_0)$, and define the corresponding piecewise-constant interpolation $t\mapsto (C^n(t),K^n(t),\boldsymbol{y}^n(t))$ as in \eqref{eq:pc1}--\eqref{eq:pc2}. \EEE
Then, we have the following:
\begin{enumerate}[(i)]
	\item For every $n\in \N$,  it holds that 
	\begin{equation*}
		\mathcal{F}(t,C^n(t),K^n(t),\boldsymbol{y}^n(t))\leq \mathcal{F}(0,C_0,K_0,\boldsymbol{y}_0)-\int_0^t \int_{\Omega} \dt{\boldsymbol{f}}(s) \cdot \boldsymbol{y}^n(s)\,\d\boldsymbol{x}\,\d s+2M\omega(\Delta_n) \quad \text{for all $t\in [0,T]$,} 
	\end{equation*}
	where $\omega\colon [0,T] \to  [0,+\infty)$ is a modulus of integrability for the function $t\mapsto \| \dt{\boldsymbol{f}}  (t)\|_{L^1(\Omega;\Rt)}$,  so that
	\begin{equation*}
		\int_S \| \dt{\boldsymbol{f}}  (s)\|_{L^1(\Omega;\Rt)}\,\d s\leq \omega(\mathscr{L}^1(S)) \quad \text{for all measurable sets $S\subset [0,T]$.}  
	\end{equation*}
	\item \EEE We have 
	\begin{equation*}
		 \sup_{t\in [0,T]} \sup_{n\in \N} \mathcal{F}(t,C^n(t),K^n(t),\boldsymbol{y}^n(t))<+\infty. 
	\end{equation*}
	\item There exists a  not relabeled   subsequence  and a function 	$t\mapsto (C(t),K(t))$ from $[0,T]$ to $\mathcal{X}_m(\closure{\Omega})$ with  both   $t\mapsto C(t)$ and $t\mapsto K(t)$ increasing such that 
	\begin{equation}
		\label{eq:helly}
		\text{$(C^{n}(t),K^{n}(t))  \to (C(t),K(t))$ in $\mathcal{X}_m(\closure{\Omega})$,   as   $n\to \infty$,     for all $t\in [0,T]$.}
	\end{equation}
	\item For all $t\in [0,T]$, there exists a  further   subsequence indexed by   $(n^t_k)_k$   and  $\boldsymbol{y}(t)\in \mathcal{Y}_{p,M}(C(t),K(t))$ for which
	\begin{equation}
		\label{eq:comp-def1}
		\text{$\boldsymbol{y}^{n^{t}_k}(t)\wks \boldsymbol{y}(t)$ in $L^\infty(\Omega;\Rt)$,} \quad \text{$D\boldsymbol{y}^{n^{t}_k}(t)\wk D\boldsymbol{y}(t)$ in $L^p(\Omega;\Rtt)$,}
	\end{equation}
	\begin{equation}
		\label{eq:comp-def2}
		\text{$\det D\boldsymbol{y}^{n^{t}_k}(t)\wk \det D\boldsymbol{y}(t)$ in $L^1(\Omega)$,}
	\end{equation}
	and 
	\begin{equation}
		\label{eq:comp-def3}
		\int_\Omega \dt{\boldsymbol{f}}(t)\cdot \boldsymbol{y}^{n_k^t}(t)\,\d\boldsymbol{x} \to \liminf_{n \to \infty} \int_\Omega \dt{\boldsymbol{f}}(t)\cdot \boldsymbol{y}^{n}(t)\,\d\boldsymbol{x},
	\end{equation}
	as $k\to \infty$. 
\end{enumerate}
\end{proposition}
\begin{proof}
First, observe that the existence of a solution $ ( (C^n_i,K^n_i,\boldsymbol{y}^n_i)  )_{i=1,\dots,\nu_n}$ of the incremental minimization problem determined by $(t^n_i)_{i=0,\dots,\nu_n}$  for every $n\in \N$  is  ensured by Proposition \ref{prop:imp}(i).

(i) Let $n\in \N$ and define $\tau^n\colon [0,T] \to \{t^n_0,\dots,t^n_{\nu_n}\}$ as in \eqref{eq:tau}.  
 By Proposition \ref{prop:imp}(ii), for all $t\in [0,T]$,  we have
\begin{equation*}
	\mathcal{F}(\tau^n(t),C^n(t),K^n(t),\boldsymbol{y}^n(t))\leq  \mathcal{F}(0,C_0,K_0,\boldsymbol{y}_0) - \int_0^{\tau^n(t)} \int_{\Omega} \dt{\boldsymbol{f}}(s)\cdot \boldsymbol{y}^n(s)\,\d\boldsymbol{x}\,\d s.
\end{equation*}
Thus, defining $\eta^n\colon [0,T]\to \R $ as
\begin{equation*}
	\eta^n(t)\coloneqq \mathcal{F}(t,C^n(t),K^n(t),\boldsymbol{y}^n(t))-\mathcal{F}(\tau^n(t),C^n(t),K^n(t),\boldsymbol{y}^n(t))+\int_{\tau^n(t)}^t \int_\Omega \dt{\boldsymbol{f}}(s)\cdot \boldsymbol{y}^n(s)\,  \d\boldsymbol{x}\,  \d s,
\end{equation*} 
we find
\begin{equation*}
	\mathcal{F}(t,C^n(t),K^n(t),\boldsymbol{y}^n(t)) + \int_0^t \int_{\Omega} \dt{\boldsymbol{f}}(\tau) \cdot \boldsymbol{y}^n(\tau)\,\d\boldsymbol{x}\,\d\tau \leq \mathcal{F}(0,C_0,K_0,\boldsymbol{y}_0) + |\eta^n(t)|.
\end{equation*}
Let $\omega$ be a modulus of integrability for the function $t\mapsto \|\dt{\boldsymbol{f}}(t)\|_{L^1(\Omega;\Rt)}$  and recall \eqref{eq:size-partition}. 
By applying the fundamental theorem of calculus for Bochner integrals,  taking advantage of the confinement condition, \EEE   we see that 
\begin{equation*}
	|\eta^n(t)|\leq 2M \int_{\tau^n(t)}^t \|\dt{\boldsymbol{f}}(s)\|_{L^1(\Omega;\Rt)}\,\d s \leq 2 M \omega(\Delta_n),
\end{equation*}
which yields the desired estimate. 

(ii) Let $n\in \N$ and $t\in [0,T]$. From (i), using the confinement condition, we obtain
\begin{equation*}
	\mathcal{F}(t,C^n(t),K^n(t),\boldsymbol{y}^n(t)) \leq  \mathcal{F}(0,C_0,K_0,\boldsymbol{y}_0)+3M \|\dt{\boldsymbol{f}}\|_{L^1(0,T;L^1(\Omega;\Rt))}.
\end{equation*}

(iii) By \cite[Theorem 6.3]{dalmaso.toader}, there exist a subsequence indexed by $(n_k)_k$ and two increasing functions $t\mapsto C(t)$ and $t\mapsto K(t)$ from $[0,T]$ to $\mathcal{K}(\closure{\Omega})$ for which \eqref{eq:helly} holds. As $\mathcal{K}_m(\closure{\Omega})$ is closed, we see that $K(t)$ 
 also lies in $\mathcal{K}_m(\closure{\Omega})$ for all $t \in [0,T]$.  From (ii), we know that 
\begin{equation*}
	\sup_{t\in [0,T]} \sup_{n\in \N} \H^0(C^n(t))<+\infty.
\end{equation*}
Thus, by Lemma \ref{lem:hausdorff-finite}, we see that $\H^0(C(t))<+\infty$, and hence $C(t)\in\mathcal{C}(\closure{\Omega})$, for all $t\in [0,T]$.

(iv) Fix $t\in [0,T]$. The claim follows from (ii) by applying Proposition \ref{prop:coercivity} to  $( (C^n(t),K^n(t),\boldsymbol{y}^n(t))   )_n$   and by standard properties of the inferior limit.  
\end{proof}

\subsection{Transfer of cavities and cracks}
In this subsection, we  prove  a variant of the so-called jump transfer lemma, see, e.g., \cite{dalmaso.francfort.toader,dalmaso.lazzaroni,francfort.giacomini.lopez}, which, 
 in the language of rate-independent systems \cite{mielke.roubicek}, corresponds to proving the existence of mutual recovery sequences.  Here, the construction needs to account also for cavitation.  

\begin{theorem}[Transfer of cavities and cracks]
	\label{thm:transfer}
Let $ ( (C_n,K_n)    )_n$ be a sequence in $  \mathcal{X}_{m}(\closure{\Omega})$    and   let  $(C,K)\in \mathcal{X}_{m}(\closure{\Omega})$ be such that $(C_n,K_n)\to (C,K)$ in $\mathcal{X}_{m}(\closure{\Omega})$. 
Also, let  $(\widetilde{C},\widetilde{K},\widetilde{\boldsymbol{y}})\in  \mathcal{A}^m_{p,M}  $ with  $\widetilde{C}\supset C$,  $\widetilde{K}\supset K$, and $ \mathcal{W}(\widetilde{\boldsymbol{y}})  <+\infty$.
Then, up to subsequences,   there exists a sequence  $  ( (\widetilde{C}_n,\widetilde{K}_n,\widetilde{\boldsymbol{y}}_n)  )_n$ in $\mathcal{A}^m_{p,M}$  satisfying the following properties:
\begin{enumerate}[(i)]
	\item $\widetilde{C}_n \supset C_n$  and  $\widetilde{C}_n\setminus C_n=\widetilde{C}\setminus C$ \EEE  for all $n\in \N$, and   $ \widetilde{C}_n   \to  \widetilde{C}  $ in   $\mathcal{K}(\closure{\Omega})$, as $n\to \infty$;
	\item $\widetilde{K}_n\supset K_n$ for all $n\in \N$,  and 
	 $\widetilde{K}_n\to \widetilde{K}$   in $\mathcal{K}(\closure{\Omega})$, and $\H^1(\widetilde{K}_n\setminus K_n)\to \H^1(\widetilde{K}\setminus K)$,  as $n\to \infty$;
	\item $\widetilde{\boldsymbol{y}}_n \to \widetilde{\boldsymbol{y}}$ and  $D\widetilde{\boldsymbol{y}}_n \to D \widetilde{\boldsymbol{y}}$ in measure  in $\Omega$, and $\mathcal{W}( \widetilde{\boldsymbol{y}}_n)\to\mathcal{W}(\widetilde{\boldsymbol{y}})$,    as $n\to \infty$;
	\item $\mathcal{P}(\widetilde{C}_n,\widetilde{K}_n,\widetilde{\boldsymbol{y}}_n)= \mathcal{P}(\widetilde{C},\widetilde{K},\widetilde{\boldsymbol{y}})$ for all $n\in \N$.
\end{enumerate}
\end{theorem}
\begin{proof}
First, we consider the case in which $\widetilde{K}$ is connected. Subsequently, we address the general case.

\textbf{\em Case 1 (Connected crack competitor).}  We subdivide the proof into  five  steps.

\emph{Step 1 (Construction of the  sets of cavitation points).} By Lemma \ref{lem:hausdorff-finite}, up to not relabeled subsequences,  there exist two finite sets $J\subset I \subset \N$  with $\H^0(C_n)=\H^0(I)$ for all $n\in \N$ and $\H^0(C)=\H^0(J)$ such that we can write  $C_n=\{\boldsymbol{a}^n_i:\: i\in I  \}$ and $C=\{ \boldsymbol{a}_j:\:j\in J \}$ in such a way that claim (ii) of Lemma \ref{lem:hausdorff-finite} holds true. 
By assumption, $\widetilde{C}\supset C$.  We define $(\widetilde{C}_n)_n$ by setting $\widetilde{C}_n\coloneqq C_n \cup \widehat{C}$ for all $n\in \N$, where $\widehat{C}\coloneqq    \widetilde{C}   \setminus C$.

\emph{Step 2 (Construction of the cracks).} We  assume that  $\widetilde{K}\in \mathcal{K}_1(\closure{\Omega})$. 
By  \cite[Lemma 3.8]{dalmaso.toader}, we find a sequence $(H_n)_n$ in $\mathcal{K}_1(\closure{\Omega})$ with $H_n\supset K_n$ for all $n\in \N$ such that
\begin{equation}
	\label{eq:H1}
	\text{$H_n\to \widetilde{K}$ in $\mathcal{K}_1(\closure{\Omega})$}
\end{equation}
and
\begin{equation}
	\label{eq:H2}
	\text{$\H^1(H_n\setminus K_n)\to \H^1(\widetilde{K}\setminus K)$. }
\end{equation}
Let $\varepsilon>0$ and set $U^\varepsilon\coloneqq \{ \boldsymbol{x}\in  \Rt  \colon \, \dist(\boldsymbol{x};\widetilde{K})<\varepsilon   \}$.  By Hausdorff convergence, there exists $\bar{n}^\varepsilon\in \N$ such that $H_n \subset  U^\varepsilon$ for $n\geq \bar{n}^\varepsilon$.  We  apply  Lemma \ref{lem:separation} with $\widetilde{K}$ in place of $K$  and with the sequence $(H_n)_n$ in place of $(K_n)_n$. Let  $N\subset \widetilde{K}$ be the corresponding set with  $\H^1( N  )=0$. Up  to  adding an $\haus$-negligible set to $N$, we additionally have
\begin{equation*}
	\lim_{r \to 0^+} \frac{\H^1(\widetilde{K} \cap Q_{\boldsymbol{\nu}_{\boldsymbol{x}}}(\boldsymbol{x},r))}{2r}=1 \quad \text{for all $\boldsymbol{x}\in \widetilde{K}\setminus N$.}
\end{equation*}
Thus, up to replacing  $\bar{r}^\varepsilon_{\boldsymbol{x}}$  in Lemma \ref{lem:separation}   with a  smaller  length,  there holds
\begin{equation}
	\label{eq:density}
 2 (1-\varepsilon)r\leq \H^1(\widetilde{K}\cap Q_{\boldsymbol{\nu}_{\boldsymbol{x}}}(\boldsymbol{x},r))\leq  2   (1+\varepsilon)r \quad \text{for all $\boldsymbol{x}\in \widetilde{K}\setminus N$ and all $r\in \big (0,\bar{r}^\varepsilon_{\boldsymbol{x}}\big)$.}
\end{equation}
We consider the family of cubes $\left\{ Q_{\boldsymbol{\nu}_{\boldsymbol{x}}}(\boldsymbol{x},r)\subset \subset U^\varepsilon:\: \boldsymbol{x}\in \widetilde{K}\setminus N,\:\: r\in \big ( 0,\bar{r}^\varepsilon_{\boldsymbol{x}}  \big)  \right\}$.   By applying the Besicovitch  covering  lemma, we find a disjoint family $\{ Q_l^\varepsilon  \}_{l= 1,\dots,  \ell^\varepsilon}$ of cubes compactly contained in $U^\varepsilon$, where $Q^\varepsilon_l\coloneqq Q_{\boldsymbol{\nu}_l^\varepsilon}(\boldsymbol{x}^\varepsilon_l,r^\varepsilon_l)$, $\boldsymbol{x}^\varepsilon_l\in \widetilde{K}\setminus N$, $\boldsymbol{\nu}_l^\varepsilon\coloneqq \boldsymbol{\nu}_{\boldsymbol{x}_l^\varepsilon}\in \So$, and $r^\varepsilon_l \in \big(0,\bar{r}^\varepsilon_{\boldsymbol{x}^\varepsilon_l}\big)$, such that
\begin{equation}
	\label{eq:besicovitch}
	\H^1\left (\widetilde{K}\setminus \bigcup_{l=1}^{\ell^\varepsilon}Q^\varepsilon_l\right )<\varepsilon.	
\end{equation}
Also, recalling \eqref{eq:H1},  by Lemma  \ref{lem:separation} we get  a sequence $(\widehat{H}^\varepsilon_n)_n$ in $\mathcal{K}_1(\closure{\Omega})$ such that
\begin{equation*}
	H_n\subset \widehat{H}^\varepsilon_n, \quad \widehat{H}^\varepsilon_n\setminus H_n \subset \bigcup_{l=1}^{\ell^\varepsilon} Q^\varepsilon_l, \quad \H^1(\widehat{H}^\varepsilon_n\setminus H_n)\leq 3\varepsilon \sum_{l=1}^{\ell^\varepsilon}r^\varepsilon_l \quad \text{for all $n\in \N$} 
\end{equation*}
and, up to replacing $\bar{n}^\varepsilon$ with a possibly larger integer and
setting $R^{\pm,\varepsilon}_l\coloneqq R^{\pm,\varepsilon}_{\boldsymbol{x}_l^\varepsilon,r^\varepsilon_l}$ according to \eqref{eq:R+}--\eqref{eq:R-},  for all $n\geq \bar{n}^\varepsilon$   we have that  
\begin{equation}
	\label{eq:sep}
	\begin{split}
		\text{$R^{+,\varepsilon}_l$ and $R^{-,\varepsilon}_l$ belong to different connected components of  $Q^\varepsilon_l \setminus \widetilde{K}$ and \EEE $Q^\varepsilon_l\setminus \widehat{H}^\varepsilon_n$ for  $l=1,\dots,\ell^\varepsilon$. }
	\end{split}
\end{equation} 
Now, let  $\lambda>1$ and  choose   $\varepsilon<\frac{\lambda-1}{5\lambda+1}$.   For $l\in \{1,\dots,\ell^\varepsilon\}$, we introduce
\begin{equation*}
	\Gamma^\varepsilon_l\coloneqq \left\{ \boldsymbol{x}\in \partial Q^\varepsilon_l:\:|(\boldsymbol{x}-\boldsymbol{x}^\varepsilon_l)\cdot \boldsymbol{\nu}^\varepsilon_l|\leq b^\varepsilon_l  \right\},
\end{equation*}
 where we set $a^\varepsilon_l\coloneqq \varepsilon r^\varepsilon_l$ and   $b^\varepsilon_l\coloneqq \frac{5\lambda+1}{\lambda-1}a^\varepsilon_l$.  Eventually, we define
 \begin{equation}
 	\label{eq:K-tilde-eps-n}
 	\widetilde{K}^\varepsilon_n\coloneqq \widehat{H}^\varepsilon_n \cup \left( \widetilde{K}\setminus \bigcup_{l=1}^{\ell^\varepsilon} Q^\varepsilon_l  \right) \cup \bigcup_{l=1}^{\ell^\varepsilon} \Gamma^\varepsilon_l.
 \end{equation} 
 We observe that $\widetilde{K}^\varepsilon_n \subset U^\varepsilon$, so that 
 \begin{equation}
 	\label{eq:hausdorff1}
 	\dist(\boldsymbol{x};\widetilde{K})<\varepsilon \quad \text{for all $\boldsymbol{x}\in \widetilde{K}^\varepsilon_n$,}
 \end{equation}
 while
 \begin{equation}
 	\label{eq:hausdorff2}
 	\dist(\boldsymbol{x};\widetilde{K}^\varepsilon_n)\leq \dist(\boldsymbol{x};H_n) \quad \text{for all $\boldsymbol{x}\in \widetilde{K}$} 
 \end{equation}
because $H_n\subset \widetilde{K}^\varepsilon_n$. Then, we estimate the length of $\widetilde{K}^\varepsilon_n\setminus K_n$ as
\begin{equation}
	\label{eq:lenght1}
	\haus(\widetilde{K}^\varepsilon_n\setminus K_n)\leq \haus(\widetilde{K}^\varepsilon_n \setminus \widehat{H}^\varepsilon_n)+\haus(\widehat{H}^\varepsilon_n \setminus H_n)+\H^1(H_n\setminus K_n).
\end{equation} 
For the first addend on the right-hand side of \eqref{eq:lenght1}, we have
\begin{equation}
	\label{eq:length2}
	\haus(\widetilde{K}^\varepsilon_n \setminus \widehat{H}^\varepsilon_n)\leq \haus\left(\widetilde{K}\setminus \bigcup_{l=1}^{\ell^\varepsilon}Q^\varepsilon_l\right)+\sum_{l=1}^{\ell^\varepsilon} \haus(\Gamma^\varepsilon_l)\leq \varepsilon + 4 \sum_{l=1}^{\ell^\varepsilon} b^\varepsilon_l \leq  \left( 1+  \frac{ 2  (5\lambda+1)}{(\lambda-1)  (1-\varepsilon)   }  \haus(\widetilde{K}) \right)\varepsilon
\end{equation} 
thanks to \eqref{eq:density}--\eqref{eq:besicovitch} and the disjointness of the cubes $\{Q^\varepsilon_l\}_l$.   Similarly, we estimate the second addend on the right-hand side of \eqref{eq:lenght1} as
\begin{equation}
	\label{eq:length3}
	\haus(\widehat{H}^\varepsilon_n \setminus H_n)\leq 3\varepsilon \sum_{l=1}^{\ell^\varepsilon} r^\varepsilon_l\leq \frac{3\varepsilon}{ 2(1-\varepsilon)} \H^1(\widetilde{K}).
\end{equation} \EEE 

\emph{Step 3 (Construction of the deformations).} We write $C_{\widetilde{\boldsymbol{y}}}\cap C=\{ \boldsymbol{a}_{j}:\:j\in J'  \}$  for some $J'\subset J$  with $  \H^0(C_{\widetilde{\boldsymbol{y}}} \cap {C})=\H^0(J')$.  We choose $\varepsilon>0$ satisfying 
\begin{equation}\label{XXX}
  2	\varepsilon  <  \dist \big(C_{\widetilde{\boldsymbol{y}}}; \partial \Omega \cup  \widetilde{K} \big) \wedge \min \big\{|\boldsymbol{b}_1-{\boldsymbol{b}}_2|: \:\boldsymbol{b}_1,{\boldsymbol{b}}_2\in C_{\widetilde{\boldsymbol{y}}},\:\boldsymbol{b}_1\neq {\boldsymbol{b}}_2 \big\}. 
\end{equation}
By  claim (ii) of Lemma \ref{lem:hausdorff-finite},   there exist a  set of indices $I'\subset I$ with $\H^0(I')=\H^0(J')$ and a bijective map $\zeta\colon I' \to J'$ such that, for all $i\in I'$, there holds $\boldsymbol{a}^n_{i}\to \boldsymbol{a}_{\zeta(i)}$, as $n\to \infty$.
Up to replacing $\bar{n}^\varepsilon$ with a larger integer, we can assume that $|\boldsymbol{a}^n_{i}-\boldsymbol{a}_{\zeta(i)}|<\varepsilon$ for all $i\in I'$ and  $n\geq \bar{n}^\varepsilon$. In this way,  by using \eqref{XXX} we find  $\{\boldsymbol{a}^n_{i}:\:i\in I'   \} \cup (C_{\widetilde{\boldsymbol{y}}}\cap \widehat{C})\subset \Omega \setminus      \closure{U^\varepsilon}$ for all $n\geq \bar{n}^\varepsilon$.  

For the construction of the deformations near the cavities, we employ the arguments in \cite[Proposition 10.1]{mora.corral}.  
Let  $n\geq \bar{n}^\varepsilon$. Consider a smooth function $f^\varepsilon \colon [0,+\infty) \to [0,1]$ satisfying
\begin{equation*}
	\text{$f^\varepsilon(z)=1$ \quad for all $ z\in [0,\varepsilon/2]$}, \qquad \text{$f^\varepsilon(z)=0$\quad for all  $ z\in  [3\varepsilon/4,+\infty)$. }
\end{equation*}
For $i\in I'$, we  define $\boldsymbol{\Phi}^\varepsilon_{n,i}\colon  \closure{B}(\boldsymbol{a}^n_{i},\varepsilon) \to \Rt$ by setting 
\begin{equation}
	\label{eq:Phi}
	\boldsymbol{\Phi}^\varepsilon_{n,i}(\boldsymbol{x})\coloneqq \boldsymbol{x}+ f^\varepsilon(|\boldsymbol{x}-\boldsymbol{a}^n_{i}|)(\boldsymbol{a}_{\zeta(i)}-\boldsymbol{a}^n_{i}). 
\end{equation}
We observe that  $\boldsymbol{\Phi}^\varepsilon_{n,i}(\boldsymbol{a}^n_{i})=\boldsymbol{a}_{\zeta(i)}$, $\boldsymbol{\Phi}^\varepsilon_{n,i}(B(\boldsymbol{a}^n_{i},\varepsilon))=B(\boldsymbol{a}^n_{i},\varepsilon)$, and also
\begin{equation}
	\label{eq:Phi-boundary}
	\boldsymbol{\Phi}^\varepsilon_{n,i}\restr{S(\boldsymbol{a}^n_{i},\varepsilon)}=\boldsymbol{id}, \qquad  (D\boldsymbol{\Phi}^\varepsilon_{n,i})\restr{S(\boldsymbol{a}^n_{i},\varepsilon)}=\boldsymbol{I}, \qquad  (D^2\boldsymbol{\Phi}^\varepsilon_{n,i})\restr{S(\boldsymbol{a}^n_{i},\varepsilon)}=\boldsymbol{O}.  
\end{equation}
 Up to replacing $\bar{n}^\varepsilon$ with a possibly larger integer,  we have that $\boldsymbol{\Phi}^\varepsilon_{n,i}$ is a diffeomorphism for $n\geq \bar{n}^\varepsilon$. Moreover
\begin{equation}
	\label{eq:Phi-determinant}
	\text{$\det D \boldsymbol{\Phi}^\varepsilon_{n,i}(\boldsymbol{x})>0$ for all $\boldsymbol{x}\in B(\boldsymbol{a}^n_{i},\varepsilon)$ and all $n\geq \bar{n}^\varepsilon$.}
\end{equation}
We define $\boldsymbol{u}^\varepsilon_{n,i} \in W^{1,p}(B(\boldsymbol{a}^n_{i},\varepsilon);\Rt)$ as $\boldsymbol{u}^\varepsilon_{n,i}\coloneqq \widetilde{\boldsymbol{y}}\circ \boldsymbol{\Phi}^\varepsilon_{n,i}$. Thus, from \cite[Proposition 5.2]{mora.corral}, we see that 
\begin{equation}
	\label{eq:Phi-cavities}
	C_{\boldsymbol{u}^\varepsilon_{n,i}}=\{ \boldsymbol{a}^n_{i}  \}, \qquad \imt(\boldsymbol{u}^\varepsilon_{n,i},\boldsymbol{a}^n_{i})=\imt(\widetilde{\boldsymbol{y}},\boldsymbol{a}_{\zeta(i)}).
\end{equation}

For the construction of deformations near the cracks, we modify the stretching argument in \cite[Lemma 4.1]{dalmaso.lazzaroni}. Let $l\in \{1,\dots,\ell^\varepsilon \}$.
For convenience,  we introduce the sets (see Figure \ref{fig:EFG}):
\begin{align*}
	F^{+,\varepsilon}_l&\coloneqq \left\{ \boldsymbol{x}\in Q^\varepsilon_l:\:-a^\varepsilon_l<(\boldsymbol{x}-\boldsymbol{x}^\varepsilon_l)\cdot \boldsymbol{\nu}_l^\varepsilon < b^\varepsilon_l \right\}, \hspace{20pt} F^{-,\varepsilon}_l\coloneqq \left\{ \boldsymbol{x}\in Q^\varepsilon_l:\:-b^\varepsilon_l<(\boldsymbol{x}-\boldsymbol{x}^\varepsilon_l)\cdot \boldsymbol{\nu}_l^\varepsilon < a^\varepsilon_l \right\},\\
	G^{+,\varepsilon}_l&\coloneqq \left\{ \boldsymbol{x}\in Q^\varepsilon_l:\:a^\varepsilon_l<(\boldsymbol{x}-\boldsymbol{x}^\varepsilon_l)\cdot \boldsymbol{\nu}_l^\varepsilon < b^\varepsilon_l \right\}, \hspace{28pt} G^{-,\varepsilon}_l\coloneqq \left\{ \boldsymbol{x}\in Q^\varepsilon_l:\:-b^\varepsilon_l<(\boldsymbol{x}-\boldsymbol{x}^\varepsilon_l)\cdot \boldsymbol{\nu}_l^\varepsilon < -a^\varepsilon_l \right\}.
\end{align*}

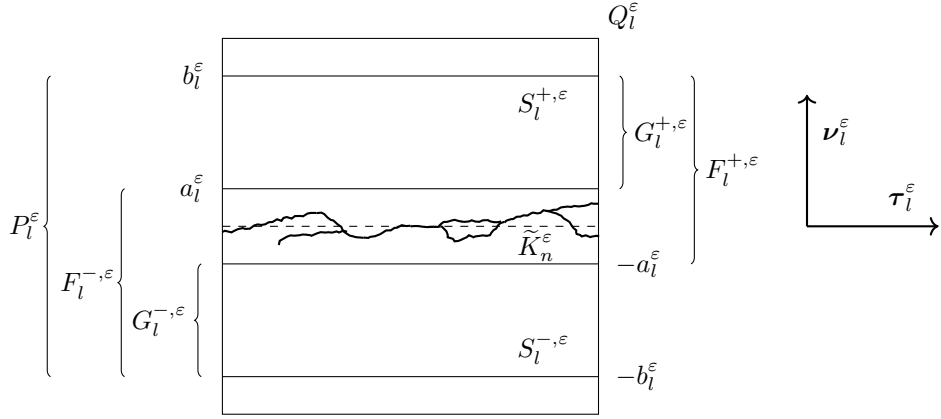
\begin{figure}
	\begin{tikzpicture}[scale=2.5]
		\draw (-1,-1) rectangle  (1,1);
		\node[above right] at (1,1) {$Q^\varepsilon_l$};
		\draw[dashed] (-1,0)--(1,0);
		\draw (-1,.2)--(1,.2);
		\draw (-1,.8)--(1,.8);
		\draw (-1,-.2)--(1,-.2);
		\draw (-1,-.8)--(1,-.8);
		\draw[decoration={brace,raise=35pt},decorate] (1,.8) -- node[right=37pt] {$F^{+,\varepsilon}_l$} (1,-.2);
		\draw[decoration={brace,raise=8pt},decorate] (1,.8) -- node[right=10pt] {$G^{+,\varepsilon}_l$} (1,.2);
		\draw[decoration={brace,raise=37pt},decorate] (-1,-.8) -- node[left=37pt] {$F^{-,\varepsilon}_l$} (-1,.2);
		\draw[decoration={brace,raise=8pt},decorate] (-1,-.8) -- node[left=10pt] {$G^{-,\varepsilon}_l$} (-1,-.2);
			\draw[decoration={brace,raise=65pt},decorate] (-1,-.8) -- node[left=65pt] {$P^{\varepsilon}_l$} (-1,.8);
		\draw[xshift=60pt,thick,->] (0,0)-- node[above right=3pt] {$\boldsymbol{\nu}^\varepsilon_l$}  (0,.7);
		\draw[xshift=60pt,thick,->] (0,0)-- node[above right=3pt] {$\boldsymbol{\tau}^\varepsilon_l$}  (.7,0);
		\node[left=3pt] at (-1,.2) {$a^\varepsilon_l$};
		\node[left=3pt] at (-1,.8) {$b^\varepsilon_l$};
		\node[right=3pt] at (1,-.2) {$-a^\varepsilon_l$};
		\node[right=3pt] at (1,-.8) {$-b^\varepsilon_l$};
		\usetikzlibrary{decorations.pathmorphing}
		\draw [ thick, decorate, decoration={random steps,segment length=1.5pt,amplitude=.6pt}]
		(-1,-.03)   to[out=0,in=180,distance=1mm] (-.5,.07);
		\draw [ thick, decorate, decoration={random steps,segment length=1.5pt,amplitude=.6pt}]
		(-.5,.07)   to[out=0,in=190,distance=1mm] (-.25,-.06);
		\draw [ thick, decorate, decoration={random steps,segment length=1.5pt,amplitude=.4pt}]
		(-.25,-.06)   to[out=0,in=170,distance=.8mm] (0,0);
		\draw [ thick, decorate, decoration={random steps,segment length=1.5pt,amplitude=.6pt}]
		(0,0)   to[out=0,in=0,distance=.4mm] (.16,0);
		\draw [ thick, decorate, decoration={random steps,segment length=1.3pt,amplitude=.8pt}]
		(.15,-.005)   to[out=350,in=200,distance=.4mm] (.24,-.07);
		\draw [ thick, decorate, decoration={random steps,segment length=1.3pt,amplitude=.6pt}]
		(.23,-.08)   to[out=0,in=130,distance=.4mm] (.33,-.062);
		\draw [ thick, decorate, decoration={random steps,segment length=1.3pt,amplitude=.4pt}]
		(.33,-.06)   to[out=10,in=90,distance=.2mm] (.4,-.06);
		\draw [ thick, decorate, decoration={random steps,segment length=1.3pt,amplitude=.3pt}]
		(.4,-.052)   to[out=10,in=180,distance=.55mm] (.51,.02);
		\draw [ thick, decorate, decoration={random steps,segment length=1.3pt,amplitude=.4pt}]
		(.5,.02)   to[out=40,in=180,distance=.7mm] (1,.12);
		\draw [ thick, decorate, decoration={random steps,segment length=1.3pt,amplitude=.4pt}]
		(.72,.083)   to[out=0,in=90,distance=.5mm] (.88,-.04);
		\draw [  thick, decorate, decoration={random steps,segment length=1.3pt,amplitude=.4pt}]
		(.88,-.04)   to[out=-20,in=10,distance=.8mm] (1,-.06);
		\draw [  thick, decorate, decoration={random steps,segment length=1.5pt,amplitude=.4pt}]
		(-.7,-.1)   to[out=75,in=20,distance=.7mm] (-0.55,-.06);
		\draw [ thick, decorate, decoration={random steps,segment length=1.4pt,amplitude=.4pt}]
		(-.53,-.052)   to[out=30,in=0,distance=.7mm] (-0.36,-.03);
		\draw [ thick, decorate, decoration={random steps,segment length=1.4pt,amplitude=.4pt}]
		(.15,0)   to[out=60,in=0,distance=.55mm] (0.25,.025);
		\draw [ thick, decorate, decoration={random steps,segment length=1.4pt,amplitude=.4pt}] 
		(.27,0.025)   to[out=10,in=-10,distance=.55mm] (0.5,.03);
		\node[right] at (.52,-.1) {$\widetilde{K}^\varepsilon_n$};
		\node[right] at (.52,.65) {${S}^{+,\varepsilon}_l$};
		\node[right] at (.52,-.65) {${S}^{-,\varepsilon}_l$};
	\end{tikzpicture}
	\caption{Representation of Step 3 of the proof of Theorem \ref{thm:transfer}. 
		The cube $Q^\varepsilon_l$ contains the sets $F^{\pm,\varepsilon}_l$, $G^{\pm,\varepsilon}_l$, and $P^\varepsilon_l$. The portion  of the cracks inside of the cube is contained in the strip $Q^\varepsilon_l \setminus (G^{+,\varepsilon}_l \cup G^{-,\varepsilon}_l)$.   Note that here  $E^{{\rm rest},\varepsilon}_{n,l}$ consists of two connected components.  }
	\label{fig:EFG}
\end{figure}

Let  again  $n\geq\bar{n}^\varepsilon$  and  define $g^\varepsilon_l\colon [-a^\varepsilon_l,b^\varepsilon_l]\to [a^\varepsilon_l,b^\varepsilon_l]$ as
\begin{equation*}
	 g^\varepsilon_l (z)\coloneqq \alpha^\varepsilon_l z^3+\beta^\varepsilon_l z^2+\gamma^\varepsilon_l z +\delta^\varepsilon_l 
\end{equation*}
by choosing  the  parameters
\begin{equation*}
	 \alpha^\varepsilon_l\coloneqq -\frac{\lambda-1}{3\lambda(a^\varepsilon_l+b^\varepsilon_l)^2}, \qquad \beta^\varepsilon_l\coloneqq \frac{(\lambda-1)b^\varepsilon_l}{\lambda(a^\varepsilon_l+b^\varepsilon_l)^2}, \qquad \gamma^\varepsilon_l\coloneqq \frac{\lambda a^\varepsilon_l(a^\varepsilon_l+2b^\varepsilon_l)+(b^\varepsilon_l)^2}{\lambda (a^\varepsilon_l+b^\varepsilon_l)^2}, \qquad \delta^\varepsilon_l\coloneqq \frac{(\lambda-1)(b^\varepsilon_l)^3}{3\lambda(a^\varepsilon_l+b^\varepsilon_l)^2}.  
\end{equation*}
The function $g^\varepsilon_l$ is smooth, strictly increasing, and convex. In particular,  recalling   $b^\varepsilon_l= \frac{5\lambda+1}{\lambda-1}a^\varepsilon_l$,  it satisfies 
\begin{equation}
	\label{eq:g-boundary}
	g^\varepsilon_l(-a^\varepsilon_l)=a^\varepsilon_l,\qquad  g^\varepsilon_l(b^\varepsilon_l)=b^\varepsilon_l, \qquad  (g^\varepsilon_l)'(-a^\varepsilon_l)=\frac{1}{\lambda},\qquad (g^\varepsilon_l)'(b^\varepsilon_l)=1, \qquad  (g^\varepsilon_l)''(b^\varepsilon_l)=0 
\end{equation}
 and
\begin{equation}
	\label{eq:g-determinant}
	\frac{1}{\lambda}\leq(g^\varepsilon_l)'(z)\leq1 \quad \text{for all $z\in [-a^\varepsilon_l,b^\varepsilon_l]$.}
\end{equation}
We consider the maps $\boldsymbol{\Psi}^{\pm,\varepsilon}_l\colon F^{\pm,\varepsilon}_l\to G^{\pm,\varepsilon}_l$  defined as
\begin{align}
	\label{eq:Psi}
	\boldsymbol{\Psi}^{\pm,\varepsilon}_l (\boldsymbol{x})\coloneqq \boldsymbol{x}^\varepsilon_l +  \left( (\boldsymbol{x}-\boldsymbol{x}^\varepsilon_l)\cdot \boldsymbol{\tau}^\varepsilon_l  \right) \boldsymbol{\tau}^\varepsilon_l \pm g^\varepsilon_l(\pm (\boldsymbol{x}-\boldsymbol{x}^\varepsilon_l)\cdot \boldsymbol{\nu}^\varepsilon_l ) \boldsymbol{\nu}^\varepsilon_l,
\end{align}
where $\boldsymbol{\tau}^\varepsilon_l\coloneqq (\boldsymbol{\nu}^\varepsilon_l \cdot \boldsymbol{e}_2,- \boldsymbol{\nu}^\varepsilon_l \cdot \boldsymbol{e}_1)^\top\in \So$ is obtained by a clockwise rotation of $\boldsymbol{\nu}^\varepsilon_l$ by $\pi/2$. We observe that $\boldsymbol{\Psi}^{\pm,\varepsilon}_l$ are smooth diffeomorphisms  for $l=1,\dots,  \ell^\varepsilon  $. In particular, from \eqref{eq:g-boundary}--\eqref{eq:g-determinant},   we have 
\begin{equation}
	\label{eq:Psi-boundary}
	\boldsymbol{\Psi}^{\pm,\varepsilon}_l\restr{S^{\pm,\varepsilon}_l}=\boldsymbol{id}, \qquad  (D\boldsymbol{\Psi}^{\pm,\varepsilon}_l)\restr{S^{\pm,\varepsilon}_l}=\boldsymbol{I}, \qquad  (D^2\boldsymbol{\Psi}^{\pm,\varepsilon}_l)\restr{S^{\pm,\varepsilon}_l}=\boldsymbol{O},
\end{equation}
where $S^{\pm,\varepsilon}_l\coloneqq \{\boldsymbol{x}\in  Q^\varepsilon_l:  \,(\boldsymbol{x}-\boldsymbol{x}^\varepsilon_l)\cdot \boldsymbol{\nu}^\varepsilon_l=\pm b^\varepsilon_l   \}$, and
\begin{equation}
	\label{eq:Psi-determinant}
	\det D \boldsymbol{\Psi}^{\pm,\varepsilon}_l(\boldsymbol{x})=  (g^\varepsilon_l)'  (\pm (\boldsymbol{x}-\boldsymbol{x}^\varepsilon_l)\cdot \boldsymbol{\nu}^\varepsilon_l )\geq \frac{1}{\lambda}\quad \text{for all $\boldsymbol{x}\in F^{\pm,\varepsilon}_l$ and $n\geq \bar{n}^\varepsilon$.}
\end{equation}
Recall \eqref{eq:sep}--\eqref{eq:K-tilde-eps-n} and set $P^\varepsilon_l\coloneqq \left\{ \boldsymbol{x}\in Q^\varepsilon_l:\,|(\boldsymbol{x}-\boldsymbol{x}^\varepsilon_l)\cdot \boldsymbol{\nu}^\varepsilon_l|<b^\varepsilon_l \right\}$. We split $P^\varepsilon_l\setminus \widetilde{K}^\varepsilon_n$ into connected components as 
\begin{equation}\label{NNN1}
	P^\varepsilon_l\setminus \widetilde{K}^\varepsilon_n=E^{+,\varepsilon}_{n,l} \cup E^{-,\varepsilon}_{n,l} \cup   E^{{\rm rest},\varepsilon}_{n,l} 
\end{equation} 
where $E^{\pm,\varepsilon}_{n,l}$ denotes the connected component containing $R^{\pm,\varepsilon}_l$ as in \eqref{eq:sep} and    $E^{{\rm rest},\varepsilon}_{n,l} $  represents the union of all other components (see Figure \ref{fig:EFG}). Then, 
we define $\boldsymbol{v}^\varepsilon_{n,l}\in W^{1,p}( (P^\varepsilon_l \cap \Omega) \EEE  \setminus \widetilde{K}^\varepsilon_n;\Rt)$ by setting 
\begin{equation}\label{NNN2}
	\boldsymbol{v}^\varepsilon_{n,l}\coloneqq \begin{cases}
		\widetilde{\boldsymbol{y}}\circ \boldsymbol{\Psi}^{+,\varepsilon}_l & \text{in $E^{+,\varepsilon}_{n,l} \cap \Omega \EEE $,}\\
		\widetilde{\boldsymbol{y}}\circ \boldsymbol{\Psi}^{-,\varepsilon}_l & \text{in $ (E^{-,\varepsilon}_{n,l}    \cup    E^{{\rm rest},\varepsilon}_{n,l} ) \cap \Omega  $.}
	\end{cases}
\end{equation}
Here, we stress that $\boldsymbol{\Psi}^{\pm,\varepsilon}_l$ takes values in $G^{\pm,\varepsilon}_l$ and $\widetilde{K} \cap Q^\varepsilon_l\subset P^\varepsilon_l \setminus (G^{+,\varepsilon}_l \cup G^{-,\varepsilon}_l)$ by \eqref{eq:sep}, \EEE   so that  $\boldsymbol{v}^\varepsilon_{n,l}$ is well defined. \EEE 

At this point, we define $\widetilde{\boldsymbol{y}}^\varepsilon_n\colon \Omega \setminus \widetilde{K}^\varepsilon_n\to \Rt$ for $n\geq\bar{n}^\varepsilon$ by setting
\begin{equation*}
	\widetilde{\boldsymbol{y}}^\varepsilon_n\coloneqq \begin{cases}
		\widetilde{\boldsymbol{y}} & \text{in $\Omega \setminus \left( \widetilde{K}^\varepsilon_n \cup \bigcup_{i\in I'} {B}(\boldsymbol{a}^n_{i},\varepsilon) \cup \bigcup_{l=1}^{\ell^\varepsilon}  {P}^\varepsilon_l  \right)$,}\\[1ex]
		\boldsymbol{u}^\varepsilon_{n,i} & \text{in $ B(\boldsymbol{a}^n_{i},\varepsilon)$ for all $i\in I'$,}\\[1ex]
		\boldsymbol{v}^\varepsilon_{n,l} & \text{in $(P^\varepsilon_l \cap \Omega )  \setminus \widetilde{K}^\varepsilon_n$ for all $l\in \{1,\dots,\ell^\varepsilon\}$.}
	\end{cases}
\end{equation*}
We need to check that  $\widetilde{\boldsymbol{y}}_n^\varepsilon \in \mathcal{Y}_{p,M}(\widetilde{C}_n,\widetilde{K}^\varepsilon_n)$  for all $n\geq \bar{n}^\varepsilon$. Clearly, $\|\widetilde{\boldsymbol{y}}_n^\varepsilon\|_{L^\infty(\Omega;\Rt)}\leq M$. The regularity of $\widetilde{\boldsymbol{y}}_n^\varepsilon$ and the positivity  of  its Jacobian determinant are verified thanks to \eqref{eq:Phi-boundary}--\eqref{eq:Phi-determinant} and \eqref{eq:Psi-boundary}--\eqref{eq:Psi-determinant}. 

In view of \eqref{eq:Phi-boundary} and \eqref{eq:Psi-boundary}, we can write $\widetilde{\boldsymbol{y}}_n^\varepsilon=\widetilde{\boldsymbol{y}}\circ \boldsymbol{d}^\varepsilon_n$ for some  injective map  $\boldsymbol{d}^\varepsilon_n \colon \Omega \setminus \widetilde{K}^\varepsilon_n \to \Rt$  of class $C^2$   with $\det D\boldsymbol{d}^\varepsilon_n(\boldsymbol{x})>0$ for all $\boldsymbol{x}\in \Omega \setminus \widetilde{K}^\varepsilon_n$. Therefore, $\widetilde{\boldsymbol{y}}_n^\varepsilon$ satisfies condition (INV) thanks to  \cite[Theorem 9.1]{mueller.spector}.  Also,  by arguing  as  in the proof of  \cite[Proposition 5.2]{mora.corral}, we show that $\mathcal{S}(\widetilde{\boldsymbol{y}}^\varepsilon_n)<+\infty$. More explicitly, given $\boldsymbol{\eta}\in C^\infty_{\rm c}((\Omega \setminus \widetilde{K}^\varepsilon_n)\times \Rt;\Rt)$ with $\|\boldsymbol{\eta}\|_{L^\infty(U\times \Rt;\Rt)}\leq 1$, we define $\widetilde{\boldsymbol{\eta}}^\varepsilon_n \in C^\infty_{\rm c}((\Omega \setminus \widetilde{K}^\varepsilon_n)\times \Rt;\Rt)$ as
\begin{equation*}
	\widetilde{\boldsymbol{\eta}}^\varepsilon_n(\boldsymbol{x},\boldsymbol{\xi})\coloneqq \begin{cases}
	\boldsymbol{\eta}(({\boldsymbol{d}}^\varepsilon_n)^{-1}(\boldsymbol{x}),\boldsymbol{\xi}) & \text{if $\boldsymbol{x} \in \boldsymbol{d}^\varepsilon_n(\pi_1(\supp \,\boldsymbol{\eta}))$,}\\
	\boldsymbol{0} & \text{if $\boldsymbol{x}\in \Omega \setminus (\widetilde{K}^\varepsilon_n \cup \boldsymbol{d}^\varepsilon_n(\pi_1(\supp\, \boldsymbol{\eta})))$,}		
	\end{cases}
\end{equation*}
where $\pi_1$ denotes the projection onto the first component. Then, applying the change-of-variable formula, we find  
\begin{equation*}
	\mathcal{S}^{\widetilde{\boldsymbol{y}}^\varepsilon_n}(\boldsymbol{\eta})=\mathcal{S}^{\widetilde{\boldsymbol{y}}}(\widetilde{\boldsymbol{\eta}}^\varepsilon_n)\leq \mathcal{S}_{\boldsymbol{d}^\varepsilon_n(\Omega \setminus \widetilde{K}^\varepsilon_n)}(\widetilde{\boldsymbol{y}})\leq \mathcal{S}(\widetilde{\boldsymbol{y}}),	
\end{equation*}
so that the claim follows.

 \EEE

By construction and by \eqref{eq:Phi-cavities}, we have
\begin{equation*}
	C_{\widetilde{\boldsymbol{y}}^\varepsilon_n}=\{\boldsymbol{a}^n_{i}:\:i\in I'   \} \cup (C_{\widetilde{\boldsymbol{y}}} \cap \widehat{C}) \subset {C}_n \cup \widehat{C}=\widetilde{C}_n
\end{equation*}
and
\begin{equation}
	\label{eq:P}
	\begin{split}
		\mathcal{P}(\widetilde{C}_n,\widetilde{K}_n,\widetilde{\boldsymbol{y}}^\varepsilon_n)&=\sum_{i\in I'} \per(\imt(\widetilde{\boldsymbol{y}}^\varepsilon_n,  \boldsymbol{a}^n_{i})  )+ \sum_{\boldsymbol{b}\in C_{\widetilde{\boldsymbol{y}}}\cap \widehat{C}} \per(\imt(\widetilde{\boldsymbol{y}},\boldsymbol{b}))\\
		&=\sum_{i\in I'} \per(\imt(\widetilde{\boldsymbol{y}},\boldsymbol{a}_{\zeta(i)})) + \sum_{\boldsymbol{b}\in C_{\widetilde{\boldsymbol{y}}}\cap \widehat{C}} \per(\imt(\widetilde{\boldsymbol{y}},\boldsymbol{b}))=\mathcal{P}(\widetilde{C},\widetilde{K},\widetilde{\boldsymbol{y}}).
	\end{split}
\end{equation}
In particular, $\widetilde{\boldsymbol{y}}_n^\varepsilon \in \mathcal{Y}_{p,M}(  \widetilde{C}_n,  \widetilde{K}^\varepsilon_n)$ as claimed.  We set
\begin{align}\label{Dneps}
	D_n^\varepsilon\coloneqq \bigcup_{i\in I'} \closure{B}(\boldsymbol{a}^n_{i},\varepsilon) \cup  \bigcup_{l=1}^{\ell^\varepsilon} \closure{P}^\varepsilon_l.
\end{align}
Clearly,
\begin{equation}
	\label{eq:D}
	\{ \widetilde{\boldsymbol{y}}^\varepsilon_n\neq \widetilde{\boldsymbol{y}}     \} \subset D^\varepsilon_n.
\end{equation} 
Moreover, 
\begin{equation}
	\label{eq:DD}
	\begin{split}
		\leb(D_n^\varepsilon)&\leq \H^0(I') \pi \varepsilon^2+4\sum_{l=1}^{\ell^\varepsilon}b^\varepsilon_lr^\varepsilon_l= \H^0(I')  \pi \varepsilon^2+ \varepsilon   \frac{4(5\lambda+1)}{\lambda-1}    \sum_{l=1}^{\ell^\varepsilon} (r^\varepsilon_l)^2\\
		&\leq \varepsilon \left(   \H^0(I') \pi \varepsilon+  \frac{4(5\lambda+1)}{\lambda-1}   \bigg(\frac{\haus(\widetilde{K})}{  2(1-\varepsilon)  }\bigg)^{\hspace{-2pt}2} \right) 	
	\end{split}
\end{equation}
thanks to \eqref{eq:density}.

\emph{Step 4 (Estimates for the elastic energy).} To estimate the  elastic energy, we employ  the  same approach  as  in  \cite[Lemma 4.1]{dalmaso.lazzaroni} and \cite[Proposition 10.1]{mora.corral}.  
 Recalling  the assumption (W4) made on  $W$, by \cite[Proposition~1.5]{dalmaso.lazzaroni}, there  exists  $\vartheta \in (0,1)$ such that
\begin{equation}
	\label{eq:W-right}
	W(\boldsymbol{F}\boldsymbol{G})\leq 2 \left( W(\boldsymbol{F})+b_W  \right) \quad \text{for all $\boldsymbol{F},\boldsymbol{G}\in \Rtt$ with $|\boldsymbol{G}-\boldsymbol{I}|<\vartheta$.}
\end{equation} 
We can choose the parameter $\lambda>1$ in Step 2  to satisfy 
\begin{equation}
	\label{eq:lambda}
	\lambda < \frac{2}{2-\vartheta}.
\end{equation}
Let $n\geq \bar{n}^\varepsilon$ and $i\in I'$. For the map $\boldsymbol{\Phi}^\varepsilon_{n,i}$ in \eqref{eq:Phi}, we have
\begin{equation*}
	D\boldsymbol{\Phi}^\varepsilon_{n,i}(\boldsymbol{x})=\boldsymbol{I}+ (f^\varepsilon)'(|\boldsymbol{x}-\boldsymbol{a}^n_{i}|) (\boldsymbol{a}_{\zeta(i)}-\boldsymbol{a}^n_{i}) \otimes \frac{\boldsymbol{x}-\boldsymbol{a}^n_{i}}{|\boldsymbol{x}-\boldsymbol{a}^n_{i}|} \quad \text{for all $\boldsymbol{x}\in B(\boldsymbol{a}^n_{i},\varepsilon)$}.
\end{equation*} 
Up to replacing $\bar{n}^\varepsilon$ with a possibly larger integer, we can assume  $|\boldsymbol{a}_{\zeta(i)}-\boldsymbol{a}^n_{i}|<{\vartheta}{ \|(f^\varepsilon)'\|^{-1}_{L^\infty(\R)}}$. In that case, we obtain
\begin{equation*}
	\|D\boldsymbol{\Phi}^\varepsilon_{n,i}-\boldsymbol{I}\|_{L^\infty(B(\boldsymbol{a}^n_{i},\varepsilon);\Rtt)}<\vartheta
\end{equation*}  
and, in turn, 
\begin{equation*}
\det D \boldsymbol{\Phi}^\varepsilon_{n,i}(\boldsymbol{x})> 1- c(\vartheta) \quad \text{for all $\boldsymbol{x}\in B(\boldsymbol{a}^n_{i},\varepsilon)$}
\end{equation*}
 for some  $0<c(\vartheta)<1$  provided that $\vartheta \ll 1$. 
Thus, using \eqref{eq:W-right},    the definition of $\boldsymbol{u}^\varepsilon_{n,i}$ before \eqref{eq:Phi-cavities},  $\boldsymbol{\Phi}^\varepsilon_{n,i}(B(\boldsymbol{a}^n_{i},\varepsilon))=B(\boldsymbol{a}^n_{i},\varepsilon)$,    and the change-of-variable formula, we estimate
\begin{equation}
	\label{eq:W-Phi}
	\begin{split}
		\int_{B(\boldsymbol{a}^n_{i},\varepsilon)} W(D\widetilde{\boldsymbol{y}}^\varepsilon_n)\,\d\boldsymbol{x}&=\int_{B(\boldsymbol{a}^n_{i},\varepsilon)} W(D\widetilde{\boldsymbol{y}}\circ \boldsymbol{\Phi}^\varepsilon_{n,i} D\boldsymbol{\Phi}^\varepsilon_{n,i})\,\d\boldsymbol{x}
		\leq \int_{B(\boldsymbol{a}^n_{i},\varepsilon)}  2  W(D\widetilde{\boldsymbol{y}}\circ \boldsymbol{\Phi}^\varepsilon_{n,i})\,\d\boldsymbol{x} +  2  b_W \pi \varepsilon^2\\
		&\leq \frac{1}{1-c(\vartheta)}\int_{B(\boldsymbol{a}^n_{i},\varepsilon)}  2  W(D\widetilde{\boldsymbol{y}}\circ \boldsymbol{\Phi}^\varepsilon_{n,i}) \det D\boldsymbol{\Phi}^\varepsilon_{n,i}\,\d\boldsymbol{x}+  2  b_W \pi \varepsilon^2\\
		&=\frac{1}{1-c(\vartheta)}\int_{B(\boldsymbol{a}^n_{i},\varepsilon)}  2   W(D\widetilde{\boldsymbol{y}})\,\d\boldsymbol{x}+ 2   b_W \pi \varepsilon^2.
	\end{split}
\end{equation}
Let $l\in \{1,\dots,\ell^\varepsilon \}$. Recalling  \eqref{eq:Psi}, we compute 
\begin{equation*}
	D \boldsymbol{\Psi}^{\pm,\varepsilon}_l(\boldsymbol{x})=\boldsymbol{\tau}^\varepsilon_l\otimes \boldsymbol{\tau}^\varepsilon_l+(g^\varepsilon_l)'((\boldsymbol{x}-\boldsymbol{x}^\varepsilon_l)\cdot \boldsymbol{\nu}^\varepsilon_l) \boldsymbol{\nu}^\varepsilon_l \otimes \boldsymbol{\nu}^\varepsilon_l=(\boldsymbol{Q}^\varepsilon_l)^\top \boldsymbol{M}(\boldsymbol{x})\boldsymbol{Q}^\varepsilon_l \quad\text{for all $\boldsymbol{x}\in F^{\pm,\varepsilon}_l$,}
\end{equation*} 
where $\boldsymbol{Q}^\varepsilon_l \in \Rtt$ is the rotation matrix having $\boldsymbol{\tau}^\varepsilon_l$ and $\boldsymbol{\nu}^\varepsilon_l$ as rows, and $\boldsymbol{M}(\boldsymbol{x})\in\Rtt$ is the diagonal matrix with entries $1$ and $(g^\varepsilon_l)'((\boldsymbol{x}-\boldsymbol{x}^\varepsilon_l)\cdot \boldsymbol{\nu}^\varepsilon_l)$. Thus, 
\begin{equation*}
	|	D \boldsymbol{\Psi}^{\pm,\varepsilon}_l(\boldsymbol{x})-\boldsymbol{I}|\leq |\boldsymbol{Q}^\varepsilon_l| \,|\boldsymbol{M}(\boldsymbol{x})-\boldsymbol{I}|\,|\boldsymbol{Q}^\varepsilon_l|\leq 2 (1-(g^\varepsilon_l)'((\boldsymbol{x}-\boldsymbol{x}^\varepsilon_l)\cdot \boldsymbol{\nu}^\varepsilon_l))\leq \frac{2(\lambda-1)}{\lambda}<\vartheta \quad\text{for all $\boldsymbol{x}\in F^{\pm,\varepsilon}_l$},
\end{equation*}
where the last two inequalities are justified by \eqref{eq:g-determinant} and \eqref{eq:lambda}, respectively. As before, we deduce
\begin{equation*}
	\det D \boldsymbol{\Psi}^{\pm,\varepsilon}(\boldsymbol{x})> 1- c(\vartheta) \quad \text{for all $\boldsymbol{x}\in F^{\pm,\varepsilon}_l$}.
\end{equation*}
 Recalling   \eqref{NNN1}--\eqref{NNN2} and   \eqref{Dneps},  with the aid of \eqref{eq:W-right} and the change-of-variable formula, we estimate
\begin{equation}
	\label{eq:W-Psi}
	\begin{split}
		\int_{P^\varepsilon_l  \cap \Omega   } W(D\widetilde{\boldsymbol{y}}^\varepsilon_n)\,\d\boldsymbol{x}&=\int_{E^{+,\varepsilon}_{n,l}  \cap \Omega  } W(D\widetilde{\boldsymbol{y}}\circ \boldsymbol{\Psi}^{+,\varepsilon}_lD\boldsymbol{\Psi}^{+,\varepsilon}_l)\,\d\boldsymbol{x} + \int_{ (E^{-,\varepsilon}_{n,l}     \cup E^{{\rm rest},\varepsilon}_{n,l}  ) \cap \Omega   } W(D\widetilde{\boldsymbol{y}}\circ \boldsymbol{\Psi}^{-,\varepsilon}_lD\boldsymbol{\Psi}^{-,\varepsilon}_l)\,\d\boldsymbol{x}\\
		&\leq  \frac{2}{1-c(\vartheta) } \int_{E^{+,\varepsilon}_{n,l}  \cap \Omega  } W(D\widetilde{\boldsymbol{y}}\circ \boldsymbol{\Psi}^{+,\varepsilon}_l)\det D \boldsymbol{\Psi}^{+,\varepsilon}_l\,\d\boldsymbol{x}\\
		 & \ \ \ + \frac{2}{1-c(\vartheta)} \int_{ ( E^{-,\varepsilon}_{n,l} \cup    E^{{\rm rest},\varepsilon}_{n,l}   )  \cap \Omega  } W(D\widetilde{\boldsymbol{y}}\circ \boldsymbol{\Psi}^{-,\varepsilon}_l)\det D \boldsymbol{\Psi}^{-,\varepsilon}_l\,\d\boldsymbol{x} +2 b_W \leb(D^\varepsilon_n)\\
		&= \frac{2}{1-c(\vartheta)} \int_{\boldsymbol{\Psi}^{+,\varepsilon}_l(E^{+,\varepsilon}_{n,l}  \cap \Omega  )\cup \boldsymbol{\Psi}^{-,\varepsilon}_l(E^{-,\varepsilon}_{n,l} \cap \Omega  )\cup   \boldsymbol{\Psi}^{-,\varepsilon}_l(    E^{{\rm rest},\varepsilon}_{n,l}    \cap \Omega   )} W(D\widetilde{\boldsymbol{y}})\,\d\boldsymbol{x} +2 b_W \leb(D^\varepsilon_n)\\
		&\leq  \frac{2}{1-c(\vartheta)} \int_{  G^{+,\varepsilon}_l \cup G^{-,\varepsilon}_l  } W(D\widetilde{\boldsymbol{y}})\,\d\boldsymbol{x}+2 b_W \leb(D^\varepsilon_n).  
	\end{split}
\end{equation}

\emph{Step 5 (Conclusion).}   The overall proof for Case 1 works as follows.   
First, as in Step 1, we define the sequence $(\widetilde{C}_n)_n$ in $\mathcal{C}(\closure{\Omega})$ by setting $\widetilde{C}_n\coloneqq C_n \cup \widehat{C}$ for every $n\in \N$, where $\widehat{C}\coloneqq \widetilde{C}\setminus C$. Then, \EEE we consider a sequence $(\varepsilon_h)_h$ of positive scalars with $\varepsilon_h \to 0^+$, as $h\to \infty$. Arguing as in Step 2 and Step 3, we find a sequence of integers $(\bar{n}^{\varepsilon_h})_h$ which we can choose to be strictly increasing. Then, we introduce another sequence $(\delta_n)_n$ by setting $\delta_n\coloneqq \varepsilon_1$ if $n\leq \bar{n}^{\varepsilon_1}$ and $\delta_n\coloneqq \varepsilon_h$ if $\bar{n}^{\varepsilon_h}\leq n < \bar{n}^{\varepsilon_{h+1}}$ for some $h\in \N$. Clearly, $\delta_n\to 0^+$, as $n\to \infty$. Following Step 2 and Step 3, we define 
the sequences $(\widetilde{K}_n)_n$,  $(\widetilde{\boldsymbol{y}}_n)_n$, and $(D_n)_n$  by setting $\widetilde{K}_n\coloneqq \widetilde{K}_n^{\delta_n}$, $\widetilde{\boldsymbol{y}}_n\coloneqq \widetilde{\boldsymbol{y}}_n^{\delta_n}$, and $D_n\coloneqq D_n^{\delta_n}$ for all $n\in \N$.

At this point, we are left to check that the sequence $ ( (\widetilde{C}_n,\widetilde{K}_n,\widetilde{\boldsymbol{y}}_n)  )_n$ in $\mathcal{A}^m_{p,M}$ satisfies all the desired properties.  By definition, $\widetilde{C}_n\setminus C_n=\widetilde{C}\setminus C$ for all $n\in\N$.  Since $C$ and $\widehat{C}$ are disjoint, and $C_n\to C$ in $\mathcal{K}(\closure{\Omega})$,  the sets $C_n$ and $\widehat{C}$ are also disjoint for $n\gg 1$.  Using the definition of the Hausdorff distance, one checks that $d(\widetilde{C}_n,\widetilde{C})\leq d(C_n,C)\to 0$, as $n\to \infty$. Thus,  (i) is proven. \EEE 
The first claim in (ii)  follows by construction;  the second one  is a consequence of   \eqref{eq:H1} and  \eqref{eq:hausdorff1}--\eqref{eq:hausdorff2}, while the third claim is established by combining \eqref{eq:H2} and \eqref{eq:lenght1}--\eqref{eq:length3}. The convergences in (iii) are consequences of   \eqref{eq:D}--\eqref{eq:DD}  and \eqref{eq:W-Phi}--\eqref{eq:W-Psi}, given that   $\mathcal{W}( \widetilde{\boldsymbol{y}})<+\infty$.   Eventually,    (iv) follows from  \eqref{eq:P}.  

\textbf{\em Case 2 (Crack competitor having multiple connected components).}   In the general case, we split $\widetilde{K}$ in connected components as $\widetilde{K}=\bigcup_{k=1}^\varkappa\widetilde{K}^k$ for  $\varkappa\in \N$ with $\varkappa\leq m$. Then, we select a   pairwise disjoint  family $\{\Omega^k\}_{k=1,\dots,\varkappa}$ of open connected subsets of $\Omega$ satisfying
\begin{equation*}
	\label{eq:V}
	\closure{\Omega}=\bigcup_{k=1}^\varkappa\closure{\Omega^k}, \qquad   \leb \left( \bigcup_{k=1}^\varkappa\partial{\Omega^k} \right)=0,  \qquad \bigcup_{k=1}^\varkappa \big( (\widetilde{C}\cap \widetilde{K}) \cap \partial \Omega^k \big) \subset \partial \Omega, 
\end{equation*}
and
\begin{equation*}
	\label{eq:VV}
  \widetilde{K} \cap \Omega^k =\widetilde{K}^k \quad \text{for all $k=1,\dots,\varkappa$.} 
\end{equation*}
Note that such a family always exists because $\H^0(\widetilde{C})<+\infty$. 
Fix $k\in \{1,\dots,\varkappa\}$. Arguing as in Case 1  for  $\Omega^k$,  $\widetilde{C}^k\coloneqq \widetilde{C}\cap \Omega^k$, and $\widetilde{K}^k$ in place of $\Omega$, $\widetilde{C}$, and $\widetilde{K}$, respectively, we find a sequence $( (\widetilde{C}^k_n, \widetilde{K}^k_n,\widetilde{\boldsymbol{y}}^k_n))_n$  for a finite set  $\widetilde{C}^k_n\subset  \closure{\Omega^k}$,  a compact set  $\widetilde{K}^k_n\subset \closure{\Omega^k}$, and $\widetilde{\boldsymbol{y}}^k_n\in W^{1,p}(\Omega^k\setminus  \widetilde{K}^k_n;  \Rt)$ for all $n\in \N$. From the result established in Case 1, we deduce that the sequence $ ( (\widetilde{C}_n,\widetilde{K}_n,\widetilde{\boldsymbol{y}}_n)  )_n$, where
\begin{equation*}
	\widetilde{C}_n\coloneqq \bigcup_{k=1}^\varkappa \widetilde{C}^k_n, \qquad \widetilde{K}_n\coloneqq \bigcup_{k=1}^\varkappa \widetilde{K}^k_n, 
\end{equation*}
and
\begin{equation*}
	\widetilde{\boldsymbol{y}}_n(\boldsymbol{x})\coloneqq \begin{cases}
		\widetilde{\boldsymbol{y}}_n^k(\boldsymbol{x}) & \text{if $\boldsymbol{x}\in \Omega^k \setminus \widetilde{K}^k_n$ for some $k \in \{1,\dots,\varkappa\},$}\\
		\widetilde{\boldsymbol{y}}(\boldsymbol{x}) & \text{if $\boldsymbol{x}\in \left( \bigcup_{k=1}^\varkappa \partial \Omega^k \right)\setminus \partial \Omega$,}
	\end{cases}
\end{equation*}
satisfies all properties (i)--(iv). In particular, for $k=1,\dots,\varkappa$, we have that $\widetilde{\boldsymbol{y}}_n^k\coloneqq \widetilde{\boldsymbol{y}} \circ \boldsymbol{d}^k_n$ for some  injective map $\boldsymbol{d}^k_n\colon \Omega^k \setminus \widetilde{K}_n \to \Rt$  of class $C^2$ with $\det D \boldsymbol{d}^k_n(\boldsymbol{x})>0$ for all $\boldsymbol{x}\in \Omega^k \setminus \widetilde{K}_n$ and $\boldsymbol{d}^k_n(\boldsymbol{x})=\boldsymbol{x}$ for all $\boldsymbol{x}  \in \partial \Omega^k$. Therefore, $\widetilde{\boldsymbol{y}}_n\in W^{1,p}(\Omega \setminus \widetilde{K}_n;\Rt)$ and $\widetilde{\boldsymbol{y}}_n$ satisfies (INV) in $\Omega \setminus \widetilde{K}_n$ by \cite[Theorem 9.1]{mueller.spector}.
\end{proof}

\subsection{Existence of irreversible quasistatic evolutions}

We are finally ready to prove our main result.

\begin{proof}[\bf Proof of Theorem \ref{thm:existence}]
We subdivide the proof into three steps.

\emph{Step 1 (Construction of the candidate solution).} Let $( (t^n_i)_{i=0,\dots,\nu_n} )_n$ be a sequence of partitions of $[0,T]$ satisfying \eqref{eq:size-partition}. For every $n\in \N$, let $  ( (C^i_n,K^i_n,\boldsymbol{y}^i_n)  )_{i=1,\dots,\nu_n}\in (\mathcal{A}^m_{p,M})^{\nu_n}$ be a solution of the incremental minimization problem determined by $(t^n_i)_{i=0,\dots,\nu_n}$ with initial datum  $(C_0,K_0,\boldsymbol{y}_0)$   whose existence is ensured by Proposition \ref{prop:imp}(i). Then, define the function $t\mapsto (C^n(t),K^n(t),\boldsymbol{y}^n(t))$ from $[0,T]$ to $\mathcal{A}^m_{p,M}$ as in \eqref{eq:pc1}--\eqref{eq:pc2}.   By Proposition \ref{prop:compactness}(iii)--(iv), there exist a function $t\mapsto (C(t),K(t),\boldsymbol{y}(t))$ from $[0,T]$ to $\mathcal{A}^m_{p,M}$ with  both  $t\mapsto C(t)$ and $t\mapsto K(t)$ increasing such that,   up to subsequences, \EEE \eqref{eq:helly} holds.  Therefore, condition (i) in Definition \ref{def:irreversible-qs} is satisfied. Moreover,  for every $t\in [0,T]$, we find a further subsequence indexed by  $(n_k^t)_k$   for which we have \eqref{eq:comp-def1}--\eqref{eq:comp-def3}. 

\emph{Step 2 (Unilateral global stability).} We prove that the function $t\mapsto (C(t),K(t),\boldsymbol{y}(t))$ satisfies condition  (ii)  in Definition \ref{def:irreversible-qs}.

 For all $n\in \N$, define the function $\tau^n\colon [0,T]\to \{t^n_0,t^n_1,\dots,t^n_{\nu_n}  \}$ as in \eqref{eq:tau}. Fix $t\in [0,T]$.  From   the minimality condition \eqref{eq:imp} satisfied by $(C^n_i,K^n_i,\boldsymbol{y}^n_i)$ for $n\in \N$ and $i\in \{1,\dots,\nu_n\}$, we deduce that  
\begin{equation}
	\label{eq:imp-n}
	\begin{split}
		&\qquad \qquad \quad \mathcal{F}(\tau^n(t),C^n(t),K^n(t),\boldsymbol{y}^n(t))\leq \mathcal{F}(\tau^n(t),\widetilde{C},\widetilde{K},\widetilde{\boldsymbol{y}})\\
		&\text{for all $n\in\N$ and for all  $(\widetilde{C},\widetilde{K},\widetilde{\boldsymbol{y}})\in\mathcal{A}^m_{p,M}$ with $\widetilde{C}\supset C^n(t)$ and $\widetilde{K}\supset K^n(t)$.}
	\end{split}
\end{equation} 
Let $(\widetilde{C},\widetilde{K},\widetilde{\boldsymbol{y}})\in\mathcal{A}^m_{p,M}$ with $\widetilde{C}\supset C(t)$ and $\widetilde{K}\supset K(t)$.  Without loss of generality, we can assume $\mathcal{W}(\widetilde{\boldsymbol{y}})<+\infty$.   
Recall \eqref{eq:helly}.  By Theorem \ref{thm:transfer}, up to subsequences, we find a sequence  $  ( (\widetilde{C}_n,\widetilde{K}_n,\widetilde{\boldsymbol{y}}_n)  )_n$ \EEE in $\mathcal{A}^m_{p,M}$ satisfying the following properties: 
\begin{enumerate}[(a)]
	\item  $\widetilde{C}_n\supset C^n(t)$   for all  $n\in \N$,  and  $\widetilde{C}_n\to \widetilde{C}$ \EEE  in $\mathcal{K}(\closure{\Omega})$  and  $\H^0(\widetilde{C}_n \setminus C^n(t))\to \H^0(\widetilde{C}\setminus C(t))$,    as $n\to \infty$;
	\item $\widetilde{K}_n\supset K^n(t)$ for all $n\in \N$, and $\widetilde{K}_n\to \widetilde{K}$ in $\mathcal{K}(\closure{\Omega})$ and   $\H^1(\widetilde{K}_n\setminus K^n(t))\to \H^1(\widetilde{K}\setminus K(t))$,   as $n\to \infty$;
	\item $\widetilde{\boldsymbol{y}}_n\to \widetilde{\boldsymbol{y}}$ and  $D\widetilde{\boldsymbol{y}}_n\to D\widetilde{\boldsymbol{y}}$ in measure  in $\Omega$, and $\mathcal{W}(\widetilde{\boldsymbol{y}}_n)\to \mathcal{W}(\widetilde{\boldsymbol{y}})$,    as $n\to \infty$;
	\item $\mathcal{P}(\widetilde{C}_n,\widetilde{K}_n,\widetilde{\boldsymbol{y}}_n)=\mathcal{P}(\widetilde{C},\widetilde{K},\widetilde{\boldsymbol{y}})$ for all $n\in \N$.
\end{enumerate}
Now, in view of (a)--(b), we can test \eqref{eq:imp-n}  with  $(\widetilde{C}_n,\widetilde{K}_n,\widetilde{\boldsymbol{y}}_n)$. \EEE We obtain  
\begin{equation*}
	 \mathcal{F}(\tau^n(t),C^n(t),K^n(t),\boldsymbol{y}^n(t))\leq \mathcal{F}( \tau^n(t),  \widetilde{C}_n, \widetilde{K}_n,\widetilde{\boldsymbol{y}}_n) 
\end{equation*}
or equivalently
\begin{equation}
	\label{eq:gus}
	\begin{split}
		 \mathcal{W}(\boldsymbol{y}^n(t))&+   \mathcal{P}(C^n(t),K^n(t),\boldsymbol{y}^n(t))-\mathcal{L}(\tau^n(t),\boldsymbol{y}^n(t))\\
		&\leq  \mathcal{W}(\widetilde{\boldsymbol{y}}_n)+\mathcal{P}(\widetilde{C}_n,\widetilde{K}_n,\widetilde{\boldsymbol{y}}_n)+\H^0(\widetilde{C}_n\setminus C^n(t))+\H^1(\widetilde{K}_n\setminus K^n(t) )-\mathcal{L}(\tau^n(t),\widetilde{\boldsymbol{y}}_n).
	\end{split}
\end{equation}
We look at the left-hand side of \eqref{eq:gus} along the sequence  in Proposition~\ref{prop:compactness}(iv).  
By applying Proposition \ref{prop:lsc}, with the aid of \eqref{eq:size-partition} and \eqref{eq:helly}--\eqref{eq:comp-def2}, we find 
\begin{equation*}
	\begin{split}
		\mathcal{W}(\boldsymbol{y}(t))+&\mathcal{P}(C(t),K(t),\boldsymbol{y}(t))-\mathcal{L}(t,\boldsymbol{y}(t))\\
		&\leq  \liminf_{k\to \infty} \left\{  \mathcal{W}(\boldsymbol{y}^{{n^t_k}}(t))+\mathcal{P}(C^{{n^t_k}}(t),K^{{n^t_k}}(t),\boldsymbol{y}^{{n^t_k}}(t)) - \mathcal{L}(\tau^{{n^t_k}}(t),\boldsymbol{y}^{{n^t_k}}(t))    \right\}.
	\end{split}
\end{equation*} 
Here, we also exploited the continuity of $\boldsymbol{f}$ from $[0,T]$ to $L^1(\Omega;\Rt)$. 
For the right-hand side of \eqref{eq:gus},   by  \eqref{eq:size-partition},  properties (a)--(d),  dominated convergence,  and, again,  by  the continuity of $\boldsymbol{f}$,   we get
\begin{equation*}
	\begin{split}
		\lim_{  n\to \infty} \big \{\EEE \mathcal{W}(\widetilde{\boldsymbol{y}}_n)+ & \mathcal{P}(\widetilde{C}_n,\widetilde{K}_n,\widetilde{\boldsymbol{y}}_n)+\H^0(\widetilde{C}_n\setminus C^n(t))+\H^1(\widetilde{K}_n\setminus K^n(t) )-\mathcal{L}(\tau^n(t),\widetilde{\boldsymbol{y}}_n)   \big \}\\
		&=\mathcal{W}(\widetilde{\boldsymbol{y}})+\mathcal{P}(\widetilde{C},\widetilde{K},\widetilde{\boldsymbol{y}})+\H^0(\widetilde{C}\setminus C(t))+\H^1(\widetilde{K}\setminus K(t) )-\mathcal{L}(t,\widetilde{\boldsymbol{y}}).
	\end{split}
\end{equation*}
Altogether, we obtain 
\begin{equation*}
	\begin{split}
		\mathcal{W}(\boldsymbol{y}(t))+&\mathcal{P}(C(t),K(t),\boldsymbol{y}(t))-\mathcal{L}(t,\boldsymbol{y}(t)) \leq \mathcal{W}(\widetilde{\boldsymbol{y}})+\mathcal{P}(\widetilde{C},\widetilde{K},\widetilde{\boldsymbol{y}})+\H^0(\widetilde{C}\setminus C(t))+\H^1(\widetilde{K}\setminus K(t) )-\mathcal{L}(t,\widetilde{\boldsymbol{y}}).
	\end{split}
\end{equation*}
Therefore, by adding $\H^0(C(t))+\H^1(K(t))$ at both sides of the previous equation,  we realize that $ (C(t),K(t),\boldsymbol{y}(t))$ satisfies the desired inequality  given in item (ii) of Definition \ref{def:irreversible-qs}.

\emph{Step 3 (Power balance).}	We show that $t\mapsto (C(t),K(t),\boldsymbol{y}(t))$ satisfies the power balance
\begin{equation*}
	\mathcal{F}(t,C(t),K(t),\boldsymbol{y}(t))= \mathcal{F}(0,C_0,K_0,\boldsymbol{y}_0)-\int_0^t \int_\Omega \dt{\boldsymbol{f}}(s)\cdot \boldsymbol{y}(s)\,\d\boldsymbol{x} \,\d s   \quad\text{for all $t\in [0,T]$.}
\end{equation*}
 From this, we  immediately   see that $t\mapsto \mathcal{F}(t,C(t),K(t),\boldsymbol{y}(t))$ belongs to $AC([0,T])$.  

Fix $t\in [0,T]$. By Proposition \ref{prop:compactness}(i), for all $n\in \N$ we have
\begin{equation}
	\label{eq:BB}
	\mathcal{F}(t,C^n(t),K^n(t),\boldsymbol{y}^n(t))\leq \mathcal{F}(0,C_0,K_0,\boldsymbol{y}_0)-\int_0^t \int_\Omega \dt{\boldsymbol{f}}(s)\cdot \boldsymbol{y}^n(s) \, \d\boldsymbol{x}\,\d s+2M\omega(\Delta_n),
\end{equation}
 where $\omega$ is a modulus of integrability of $s\mapsto \|\dt{\boldsymbol{f}}(s)\|_{L^1(\Omega;\Rt)}$.  
 Recall  Proposition \ref{prop:compactness}(iv). 
For the left-hand side of \eqref{eq:BB},   with the aid of \eqref{eq:helly}--\eqref{eq:comp-def2} and Proposition \ref{prop:lsc}, we get
\begin{equation*}
	\label{eqn:LHS}
	\mathcal{F}(t,C(t),K(t),\boldsymbol{y}(t))\leq \liminf_{k\to \infty} \mathcal{F}(t,C^{n^t_k}(t),K^{n^t_k}(t),\boldsymbol{y}^{n^t_k}(t))\leq \limsup_{n \to \infty} \mathcal{F}(t,C^n(t),K^n(t),\boldsymbol{y}^n(t)).
\end{equation*}
For the right-hand side of \eqref{eq:BB}, applying the reverse Fatou lemma, we find
\begin{equation*}
	\begin{split}
		\limsup_{n \to \infty} \int_0^t \left( - \int_\Omega \dt{\boldsymbol{f}}(s)\cdot \boldsymbol{y}^n(s)\,\d\boldsymbol{x} \right)\,\d s&\leq \int_0^t \limsup_{n \to \infty} \left( - \int_\Omega \dt{\boldsymbol{f}}(s)\cdot \boldsymbol{y}^n(s)\,\d\boldsymbol{x} \right) \,\d s\\
		&=-\int_0^t \liminf_{n \to \infty} \left(  \int_\Omega \dt{\boldsymbol{f}}(s)\cdot \boldsymbol{y}^n(s)\,\d\boldsymbol{x} \right)\,\d s\\
		&=-\int_0^t \lim_{k\to \infty} \left( \int_\Omega \dt{\boldsymbol{f}}(s)\cdot \boldsymbol{y}^{n^s_k}(s)\,\d\boldsymbol{x}\right)\,\d s\\
		&=-\int_0^t \int_\Omega \dt{\boldsymbol{f}}(s)\cdot \boldsymbol{y}(s)\,\d \boldsymbol{x}\,\d s, 
	\end{split}
\end{equation*}
where in the last two lines we resorted to \eqref{eq:comp-def1} and \eqref{eq:comp-def3}.
Hence, recalling \eqref{eq:size-partition}, the combination of the previous two inequalities yields \EEE 
\begin{equation*}
	\mathcal{F}(t,C(t),K(t),\boldsymbol{y}(t))\leq \mathcal{F}(0,C_0,K_0,\boldsymbol{y}_0)-\int_0^t \int_\Omega \dt{\boldsymbol{f}}(s)\cdot \boldsymbol{y}(s)\,\d\boldsymbol{x} \,\d s .
\end{equation*} 

To prove the opposite inequality, we adopt the  by now  classical argument based on the global stability and  the approximation of the Lebesgue integral by Riemann sums. By applying \cite[Lemma 4.12 and Remark~4.13]{dalmaso.francfort.toader} to the functions $s\mapsto \int_\Omega \dt{\boldsymbol{f}}(s)\cdot \boldsymbol{y}(s)\,\d\boldsymbol{x}$ and $s\mapsto \dt{\boldsymbol{f}}(s)$ from $[0,t]$ to $\R$ and $L^1(\Omega;\Rt)$, respectively, we find  a sequence $( (s^h_i)_{i=0,\dots,\mu_h} )_h$ of partitions of $[0,t]$ 
 with  $  \max_{ i=1,\dots, \mu_h} |s_i^h-s_{i-1}^h| \to 0 $, as $h\to \infty$,  satisfying the following: 
\begin{equation}
	\label{eq:p2}\sum_{i=1}^{\mu_h} \left | (s^h_i-s^h_{i-1}) \int_\Omega \dt{\boldsymbol{f}}(s^h_i)\cdot \boldsymbol{y}(s^h_i)\,\d\boldsymbol{x}-\int_{s^h_{i-1}}^{s^h_i} \int_\Omega \dt{\boldsymbol{f}}(s)\cdot \boldsymbol{y}(s)\,\d\boldsymbol{x}\,\d s \right |\to 0,\quad \text{as $h\to \infty$}
\end{equation}
and
\begin{align}
	\label{eq:p1}\sum_{i=1}^{\mu_h} \left \| (s^h_i-s^h_{i-1}) \dt{\boldsymbol{f}}(s^h_i)-\int_{s^h_{i-1}}^{s^h_i} \dt{\boldsymbol{f}}(s)\,\d s\right \|_{L^1(\Omega;\Rt)}\to 0, \quad \text{as $h\to \infty$.}
\end{align}
Let $h\in \N$ and $i\in \{ 1,\dots,\mu_h \}$. Given that $C(s^h_i)\supset C(s^h_{i-1})$ and $K(s^h_i)\supset C(s^h_{i-1})$, we can test the global stability of $(C(s^h_{i-1}),K(s^h_{i-1}),\boldsymbol{y}(s^h_{i-1}))$ by choosing $(C(s^h_i),K(s^h_i),\boldsymbol{y}(s^h_i))$ as a competitor. We obtain
\begin{equation*}
	\begin{split}
		&\mathcal{F}(s^h_{i-1}, C(s^h_{i-1}),K(s^h_{i-1}),\boldsymbol{y}(s^h_{i-1}))\leq \mathcal{F}(s^h_{i-1},C(s^h_{i}),K(s^h_{i}),\boldsymbol{y}(s^h_{i}))\\
		&\qquad =\mathcal{W}(\boldsymbol{y}(s^h_i))+\mathcal{P}(C(s^h_i),K(s^h_i),\boldsymbol{y}(s^h_i))+\H^0(C(s^h_i))+\H^1(K(s^h_i))-\mathcal{L}(s^h_{i-1},\boldsymbol{y}(s^h_i))\\
		&\qquad=\mathcal{F}(s^h_i,C(s^h_i),K(s^h_i),\boldsymbol{y}(s^h_i))+\int_{s^h_{i-1}}^{s^h_i} \int_\Omega \dt{\boldsymbol{f}}(s)\cdot \boldsymbol{y}(s^h_i)\,\d\boldsymbol{x}\,\d s,
	\end{split}
\end{equation*}
where in the last line we applied the fundamental theorem of calculus for Bochner integrals. Summing the previous inequality for $i=1,\dots,\mu_h$ yields
\begin{equation}
	\label{eq:uei}
	\mathcal{F}(0,C_0,K_0,\boldsymbol{y}_0)\leq \mathcal{F}(t,C(t),K(t),\boldsymbol{y}(t)) + \sum_{i=1}^{\mu_h} \int_{s^h_{i-1}}^{s^h_i} \int_\Omega \dt{\boldsymbol{f}}(s)\cdot \boldsymbol{y}(s^h_i)\,\d\boldsymbol{x}\,\d s.
\end{equation}
We rewrite the last term as
\begin{equation*}
	\begin{split}
		\sum_{i=1}^{\mu_h} \int_{s^h_{i-1}}^{s^h_i} \int_\Omega \dt{\boldsymbol{f}}(s)\cdot \boldsymbol{y}(s^h_i)\,\d\boldsymbol{x}\,\d s&= 	\sum_{i=1}^{\mu_h} (s^h_i-s^h_{i-1}) \int_\Omega \dt{\boldsymbol{f}}(s^h_i)\cdot \boldsymbol{y}(s^h_i)\,\d\boldsymbol{x}\\
		& \quad + 	\sum_{i=1}^{\mu_h} \int_{s^h_{i-1}}^{s^h_i}  \int_\Omega \left( \dt{\boldsymbol{f}}(s)-\dt{\boldsymbol{f}}(s^h_i)  \right)\cdot \boldsymbol{y}(s^h_i)\,\d\boldsymbol{x}\,\d s.
	\end{split}
\end{equation*}
As $h\to \infty$, the first term on the right-hand side of the previous equation converges to  $\int_0^t \int_\Omega \dt{\boldsymbol{f}}(s)\cdot \boldsymbol{y}(s)\,\d\boldsymbol{x}\,\d s$ by \eqref{eq:p2}, while the second one vanishes as a consequence of \eqref{eq:p1},  the confinement condition,  and standard properties of the Bochner integral. Thus, letting $h\to \infty$ in \eqref{eq:uei}, we obtain
\begin{equation*}
	\mathcal{F}(0,C_0,K_0,\boldsymbol{y}_0)-\int_0^t \int_\Omega \dt{\boldsymbol{f}}(s)\cdot \boldsymbol{y}(s)\,\d\boldsymbol{x} \leq 	\mathcal{F}(t,C(t),K(t),\boldsymbol{y}(t)),
\end{equation*}
which concludes the proof.
\end{proof}
  
\section*{Acknowledgements}
 M.~Bresciani acknowledges the support of the Alexander von Humboldt Foundation.   M.~Friedrich has been  supported by the DFG project FR 4083/5-1.

\end{document}